\documentclass[a4paper,12pt]{article}
\usepackage{dsfont}
\usepackage{amsmath}
\usepackage{enumitem}
\usepackage{amssymb,amsthm}
\usepackage{mathrsfs,amsfonts}
\usepackage{color}
\usepackage{graphicx}
\usepackage{bbm}
\usepackage[normalem]{ulem}
\usepackage{hyperref}
\usepackage[capitalise]{cleveref}

\allowdisplaybreaks

\allowdisplaybreaks

\newtheorem{con0}{Theorem}[section]
\newtheorem{thm0}{Theorem}[section]
\newtheorem{exa0}{Theorem}[section]

\newtheorem{con1}[con0]{Condition}
\newtheorem{def1}[thm0]{Definition}
\newtheorem{lem1}[thm0]{Lemma}
\newtheorem{thm1}[thm0]{Theorem}
\newtheorem{cor1}[thm0]{Corollary}
\newtheorem{pro1}[thm0]{Proposition}
\newtheorem{rem1}[thm0]{Remark}
\newtheorem{ass1}[thm0]{Assumption}
\newtheorem{exa1}[exa0]{\it{Example}}

\def\bglemma{\begin{lem1}}\def\edlemma{\end{lem1}}
\def\bgtheorem{\begin{thm1}}\def\edtheorem{\end{thm1}}
\def\btheorem{\begin{thm1}}\def\etheorem{\end{thm1}}

\def\bgproposition{\begin{pro1}}\def\edproposition{\end{pro1}}

\def\benumerate{\begin{enumerate}}\def\eenumerate{\end{enumerate}}
\def\bitemize{\begin{itemize}}\def\eitemize{\end{itemize}}

\def\beqlb{\begin{eqnarray}}\def\eeqlb{\end{eqnarray}}
\def\beqnn{\begin{eqnarray*}}\def\eeqnn{\end{eqnarray*}}

\def\eqref#1{{\rm(\ref{#1})}}

\def\ar{\!\!\!&}\def\nnm{\nonumber}\def\ccr{\nnm\\}

\def\proof{\noindent{\it Proof.~}}
\def\qed{\hfill$\square$\smallskip}

\def\mrm{\mathrm}\def\mbb{\mathbf}\def\mcr{\mathscr}
\def\mbb{\mathbf}\def\mds{\mathds}

\def\d{\mrm{d}}\def\e{\mrm{e}}

\def\I{\mds{1}}

\begin{document}

\bigskip\bigskip

\centerline{\Large\bf A Pathwise Approach to the Strong Feller Property}
\centerline{\Large\bf and Irreducibility of Nonlinear Branching Processes}

\bigskip\bigskip

\begin{center}
{\large Pei-Sen Li}\\
School of Mathematics and Statistics, Beijing Institute of Technology, Beijing 100872, China\\
\href{mailto:peisenli@bit.edu.cn}{\texttt{peisenli@bit.edu.cn}}\\[4pt]
{\large Xiangqi Zheng}\\
School of Mathematics, East China University of Science and Technology, Shanghai 201206, China\\
\href{mailto:zhengxq@ecust.edu.cn}{\texttt{zhengxq@ecust.edu.cn}}\\[4pt]
{\large Xiaowen Zhou}\\
Department of Mathematics and Statistics, Concordia University,\\
1455 De Maisonneuve Blvd. W., Montreal, Canada\\
\href{mailto:xiaowen.zhou@concordia.ca}{\texttt{xiaowen.zhou@concordia.ca}}
\end{center}

\bigskip\bigskip

\medskip

\bigskip

\begin{abstract}
We study the strong Feller property and irreducibility for continuous-state nonlinear branching processes defined as solutions to stochastic differential equations with jumps. Due to boundary degeneracy and discontinuous jump coefficients, classical methods do not apply. We develop a pathwise approach combining state-dependent time change, truncated auxiliary processes, and localized coupling to establish these two properties. As applications, we obtain exponential convergence to a unique quasi-stationary distribution in the absorbing case, and uniform exponential ergodicity in the non-absorbing case. This pathwise approach  is flexible and can be adapted to a broader class of jump-diffusions without relying on specific coefficient structures.

\smallskip
\noindent \textbf{Keywords:} Continuous-state branching processes,
degenerate jump-diffusions, strong Feller property, irreducibility, state-dependent time-change,
quasi-stationary distributions, uniform exponential ergodicity.
\end{abstract}

\bigskip

\smallskip

\noindent{\textit{MSC (2020) Subject Classification:}} 60J80, 60H10, 60J76

\bigskip

\section{Introduction}

 \setcounter{equation}{0}
Continuous-state branching processes (CB-processes) arise as scaling limits of
classical Galton--Watson branching processes \cite{Feller,Lamperti} and can be
represented via stochastic differential equations \cite{DawsonLi}.
Although classical CB-processes are characterized by the branching property,
meaning that individuals reproduce and die independently of one another,
realistic biological and ecological systems often exhibit \textit{density-dependent}
interactions when the population size is large or resources are limited.  To
capture such interactions, various generalizations have been introduced;
see, e.g., \cite{Lam05, BFF18, PW15, Le2013}.

Recently, Li, Yang, and Zhou \cite{LYZ} proposed a general class of
\textit{continuous-state nonlinear branching processes} (CSNBPs), which unifies
many of these models by allowing the branching rates, the diffusion coefficient,
and the drift to depend nonlinearly on the current population size.  On a
filtered probability space
\((\Omega,\mathcal F,(\mathcal F_t)_{t\ge0},\mathbf P)\) satisfying the usual
conditions, the process $X=(X_t)_{t\ge0}$ is defined as the unique nonnegative
strong solution to the SDE
\begin{equation}\label{eq:sde_intro}
\begin{aligned}
X_t = X_0 &+ \int_0^t R_0(X_s)\,\d s
          + \int_0^t \sqrt{R_1(X_s)}\,\d B_s   \\
          &+ \int_0^t \int_0^\infty \int_0^{R_2(X_{s-})}
             z\,\tilde N(\d s,\d z,\d u),
\end{aligned}
\end{equation}
where $(B_t)_{t\ge0}$ is a standard Brownian motion,
$N(\d s,\d z,\d u)$ is an independent Poisson random measure on
$(0,\infty)^3$ with intensity $\d s\,\pi(\d z)\,\d u$ satisfying
$\int_0^\infty (z\wedge z^2)\,\pi(\d z)<\infty$, and $\tilde N$ is the
corresponding compensated measure.  The coefficients $R_0,R_1,R_2$ are
continuous on $[0,\infty)$ and satisfy the natural local Lipschitz and boundary
conditions; the precise assumptions are stated in Section~2.

The strong Feller property guarantees that for any fixed $t>0$, the transition semigroup maps bounded measurable functions to continuous ones, thereby smoothing out the dependence on the initial state; irreducibility ensures that the process can reach any open set from any interior starting point. Together, they form the backbone of uniqueness proofs for invariant and quasi-stationary distributions and are also used in proving exponential ergodicity.

The analysis of the strong Feller property and irreducibility for the CSNBP \eqref{eq:sde_intro} is, however, complicated by three structural features of the model. (i)~\textbf{Boundary degeneracy}: both the diffusion coefficient $R_1(x)$ and the jump intensity $R_2(x)$ naturally degenerate at the boundary $x=0$. (ii)~\textbf{Discontinuous jump mechanism}: the Poisson integral contains the indicator $\mathds{1}_{\{u\le R_2(X_{s-})\}}$, which is discontinuous in the state variable. (iii)~\textbf{Absence of branching structure}: the nonlinear dependencies in $R_0$, $R_1$, and $R_2$ prevent the use of the affine structure and the classical branching property.

There is a large literature on the strong Feller property and irreducibility
for jump-diffusions.  Most general methods require the driving noise to be
nondegenerate and the coefficients to be sufficiently regular in the state
variable.  For the strong Feller property, analytic approaches such as
Bismut-type formulae \cite{PZ95,DPEZ95}, Krylov estimates, and Zvonkin's
transforms \cite{XZ20} require uniform ellipticity or sufficiently rich
L\'evy noise.  Classical global coupling methods
\cite{Q14,XZhu19,KZ20} and coupling-by-change-of-measure arguments leading to
Harnack inequalities \cite{WangYuan11a} make similar nondegeneracy and
regularity assumptions.  For irreducibility, Brownian-based arguments, such as
Girsanov-type approaches, again rely on uniform ellipticity and regular jump
coefficients \cite{PZ95,Q14,KZ20,XZhu19}.  Other approaches use L\'evy
noise through support or reachability arguments, but they are usually developed
either in additive-noise settings or under nondegenerate pure-jump structures
\cite{PSXZ12,FHR16,WYZZ22}.  These assumptions are incompatible with Features
(i) and (ii).

A separate class of results treats affine or linear branching processes with
degenerate coefficients at the boundary.  The proofs often exploit special
structural features, such as affine transform formulae, branching-immigration
decompositions, Lamperti-type time changes, or stable branching noise
structures \cite{Li2019Polynomial,LM15,lwz,ChenLi23,FJKR19,FriesenJin2020JCIR}. Feature (iii)
prevents us from using these specialized branching arguments in the nonlinear
setting considered here.

We overcome these obstacles by developing a pathwise framework that does not rely on any affine structure. Our methodological contributions are twofold.

First, we prove irreducibility by reducing the problem to one-sided reachability.
Upward movement is supplied either by the interior diffusion or by positive jumps.
For downward movement, which is delicate in the absence of negative jumps, we use a
state-dependent time change together with truncated auxiliary processes. Under the
infinite-variation small-jump condition, the truncation leaves a compensating drift
strong enough to drive the auxiliary process downward. This gives downward
reachability; together with upward reachability, it proves irreducibility on the
interior state space.

Second, we establish the strong Feller property by a localized coupling argument.
Classical global coupling estimates are not available because the noise degenerates
near the boundary and the Poisson integral contains a discontinuous state-dependent
indicator. Instead, we localize the coupling near a fixed interior point and compare
the coupling time with the exit time from a small neighbourhood. The local contraction
is produced either by the non-degenerate interior diffusion or by the high activity of
small jumps. An exit estimate then shows that two processes starting from nearby
interior points couple before leaving the neighbourhood with probability tending to one.
This yields the interior strong Feller property.

With the strong Feller property and irreducibility established under mild noise conditions (Theorems~\ref{SF_a} and~\ref{irred}), standard Lyapunov drift criteria can be applied to the long-time analysis. Under appropriate Lyapunov conditions, the behaviour splits into two boundary regimes. When the boundary is absorbing, there exists a unique quasi-stationary distribution to which the process converges exponentially fast (Theorem~\ref{qsd}); when the boundary is non-absorbing, the process is uniformly exponentially ergodic (Theorem~\ref{thmergo}). The precise statements and assumptions are given in Section~2.

\textbf{Organization of the paper.} Section~2 introduces the analytical assumptions and states the main results. Section~3 proves the strong Feller and trajectory Feller properties via localized coupling arguments. Section~4 establishes global irreducibility using the time-change technique and auxiliary paths. Finally, Section~5 applies these topological results to the QSD in the absorbing case and to uniform exponential ergodicity in the non-absorbing regime.

\section{Assumptions and Main Results}
\subsection{Standing Assumptions}

Although our primary motivation comes from the population dynamics of the CSNBP \eqref{eq:sde_intro}, the regularity and ergodicity properties we investigate inherently depend on the analytical properties of the coefficients. To highlight the general applicability of our pathwise framework to degenerate jump-diffusions, we formulate our core assumptions directly as analytical constraints.

Throughout this paper, we assume that the coefficients $R_0, R_1$, and $R_2$ are continuous functions on $[0, \infty)$ satisfying the following basic conditions:
\begin{itemize}
    \item For any compact interval $A\subset (0,\infty)$, there exists a constant $C(A)>0$ such that for all $x,y\in A$,
    \beqlb\label{lclip}
    |R_0(x)-R_0(y)|+|R_1(x)-R_1(y)|+|R_2(x)-R_2(y)|\le C(A)|x-y|.
    \eeqlb
    \item $R_0: [0, \infty) \to \mathbb{R}$ satisfies the boundary condition $R_0(0) \ge 0$.
    \item $R_1: [0, \infty) \to [0, \infty)$ satisfies the boundary condition $R_1(0) = 0$.
    \item $R_2: [0, \infty) \to [0, \infty)$ is non-decreasing, strictly positive on $(0,\infty)$, and satisfies $R_2(0) = 0$.
\end{itemize}

Since our primary focus is on the topological properties of the transition semigroup, we shall not discuss the criteria for the existence of a strong solution. Instead, we operate under the standing assumption that for any initial value $X_0 \ge 0$, the SDE \eqref{eq:sde_intro} admits a pathwise unique and non-explosive strong solution taking values in $[0, \infty)$ (see \cite{LYZ} for such well-posedness criteria).

\subsection{Main results}

We now present the main theorems concerning the topological properties of the continuous-state nonlinear branching process. As highlighted in the introduction, establishing the strong Feller property and irreducibility requires stochastic fluctuations that smooth the dependence on initial conditions and ensure reachability across the state space. In our framework, these fluctuations can originate either from the continuous diffusion component or from the accumulation of active small jumps:
\begin{description}
    \item[\hypertarget{cond:C1}{\textbf{(C1)}}] $R_1$ is strictly positive on $(0,\infty)$.
    
    \item[\hypertarget{cond:C2}{\textbf{(C2)}}] There exists $\gamma>0$ and $\alpha\in(0, 2)$ such that $\pi(\d z)\ge\I_{\{0<z\le 1\}} \gamma z^{-1-\alpha}\d z$.
    
    \item[\hypertarget{cond:C3}{\textbf{(C3)}}] $\int_0^1 z\pi(\d z)=\infty$.
    \end{description}

Let $(X_t)_{t\ge0}$ be the solution of \eqref{eq:sde_intro}, and let $(P_t)_{t\ge0}$ be its transition semigroup. For probability measures \(\mu,\nu\) on \([0,\infty)\), write
\[
\|\mu-\nu\|_{\mathrm{TV}}=\sup_{B\in\mathcal{B}([0,\infty))}|\mu(B)-\nu(B)|.
\]
Our first main result establishes the strong Feller property in the interior state space.

\bgtheorem[Strong Feller property]\label{SF_a}
Suppose that either \hyperlink{cond:C1}{\textup{(C1)}} or \hyperlink{cond:C2}{\textup{(C2)}} holds. Then for any \(t>0\) and any \(x_0\in(0,\infty)\),
\[
\lim_{x\to x_0}\|P_t(x,\cdot)-P_t(x_0,\cdot)\|_{\mathrm{TV}}=0.
\]
Consequently, for any \(t>0\) and any \(f\in b\mcr B([0,\infty))\), the mapping \(x\mapsto P_tf(x)\) is continuous on \((0,\infty)\).
\edtheorem

Our second main result establishes the irreducibility of the process.

\btheorem[Irreducibility]\label{irred}
Suppose that either \hyperlink{cond:C1}{\textup{(C1)}} or \hyperlink{cond:C3}{\textup{(C3)}} holds. Then $(P_t)_{t\ge0}$ is irreducible; that is, for any $x>0$, $t>0$, and any nonempty open set $O\subset(0,\infty)$, we have $P_t(x,O)>0.$
\etheorem

\begin{rem1}\label{rem_conditions}
Conditions \hyperlink{cond:C1}{\textup{(C1)}}--\hyperlink{cond:C3}{\textup{(C3)}} should be read as replacements for the global nondegeneracy and regularity assumptions commonly imposed in the literature.

First, Condition \hyperlink{cond:C1}{\textup{(C1)}} is strictly weaker than uniform ellipticity: it allows the diffusion coefficient to degenerate at the boundary $x=0$, while retaining non-degeneracy in the interior. This interior diffusion is sufficient both for local smoothing and for bidirectional movement before the process approaches the boundary.

Second, when smoothing is supplied by the jump noise rather than by the diffusion, Condition \hyperlink{cond:C2}{\textup{(C2)}} provides a stable-like lower bound on the small-jump activity. This replaces the usual smoothness or Lipschitz assumptions on the jump coefficient, appearing for instance in \cite{KZ20,KF99}, and is essential here, since the branching jump mechanism contains the discontinuous indicator $\mds{1}_{\{u\le R_2(X_{s-})\}}$.

We also note the relation between \hyperlink{cond:C2}{\textup{(C2)}} and \hyperlink{cond:C3}{\textup{(C3)}}. If the exponent in \hyperlink{cond:C2}{\textup{(C2)}} can be chosen with \(\alpha\ge1\), then \hyperlink{cond:C2}{\textup{(C2)}} implies \hyperlink{cond:C3}{\textup{(C3)}}; for \(\alpha\in(0,1)\), however, the stable-like lower bound in \hyperlink{cond:C2}{\textup{(C2)}} alone does not provide the infinite-variation input required in \hyperlink{cond:C3}{\textup{(C3)}}.

Third, for irreducibility in the absence of continuous diffusion, Condition \hyperlink{cond:C3}{\textup{(C3)}} supplies the infinite-variation input needed to compensate for the absence of negative jumps. In particular, Theorem~\ref{irred} establishes irreducibility on $(0,\infty)$ even though the jump intensity $R_2(x)$ is state-dependent and vanishes at the boundary, in contrast with approaches based on non-degenerate additive L\'evy noise; see, e.g., \cite{PSXZ12,FHR16}.
\end{rem1}

\subsection{Applications: quasi-stationarity and ergodicity}
In this section, we apply the established strong Feller property and irreducibility to investigate the long-time behaviour of CSNBPs. Under suitable Lyapunov conditions, the interplay of these properties leads to two asymptotic regimes: in the absorbing case where $R_0(0)=0$, we establish the existence and uniqueness of a quasi-stationary distribution, whereas in the non-absorbing case where $R_0(0)>0$, we obtain uniform exponential ergodicity.

Let $\mathcal{P}([0,\infty))$ denote the space of Borel probability measures on $[0,\infty)$. For $\nu\in\mathcal P([0,\infty))$, we write $\nu P_t$ for the distribution of $X_t$ with initial distribution $\nu$.
We further impose the following Lyapunov-type drift condition:

\medskip

\noindent\hypertarget{cond:C4}{{\bf(C4)}} For some $r>0$, $\lambda_1>0$, and $\lambda_2\ge0$, it holds that
\[
R_0(x)\le -\lambda_1 x^{r+1}+\lambda_2, \qquad x>0.
\]
\noindent Condition \hyperlink{cond:C4}{\textup{(C4)}} serves as a strong dissipation assumption. By requiring the drift $R_0(x)$ to exhibit super-linear negative growth for large $x$, it provides the Lyapunov drift needed to control large-state excursions. While Condition \hyperlink{cond:C4}{\textup{(C4)}} governs the boundary behavior at $\infty$, we also require an assumption to characterize the local dynamics near $0$.

\medskip

\noindent\hypertarget{cond:C5}{\textbf{(C5)}} There exist constants $k>0$, $\delta\in(0,1)$, and $\theta\in(0,1)$ such that \(R_0(0)=0\), \(R_0(x)\le kx^\theta\) for all $x\in(0,\delta)$, and at least one of the following holds on the same interval:
\begin{enumerate}[label=\rm(\roman*), leftmargin=2em, itemsep=4pt]
    \item $R_1(x)\ge \frac{2k}{\theta} x^{\theta+1}$ for all $x\in(0,\delta);$
    \item There exist constants $\gamma>0$ and $\alpha\in(0, 2)$ such that 
    \[
    \pi(\d z)\ge \gamma z^{-1-\alpha}\I_{\{0<z\le 1\}}\,\d z \quad \text{and} \quad R_2(x) \ge \frac{12(2-\alpha)}{\gamma}kx^{\theta+\alpha-1}\quad\text{for all }x\in(0,\delta).
    \]
\end{enumerate}

\noindent Condition \hyperlink{cond:C5}{\textup{(C5)}} controls the local dynamics near the lower boundary $0$. Near the origin, the drift term $R_0$ might push the process towards the interior $(0, \infty)$, away from the boundary. Condition \hyperlink{cond:C5}{\textup{(C5)}} requires that this outward drift be dominated either by a sufficiently strong continuous diffusion $R_1$ (case (i)), or by a high intensity of small jumps characterized by $\alpha$ and $R_2$ (case (ii)). This boundary condition is the key local input used later in the quasi-stationary analysis.

In population dynamics,  when the state $0$ is absorbing ($R_0(0)=0$), the standard invariant measure is merely a Dirac mass at zero, which provides no insight into the persistent behavior of the population prior to extinction. To capture the conditional long-term dynamics, we investigate the quasi-stationary distribution (QSD). Let $\tau_0^-=\inf\{t\ge0:X_t\le0\}$ be the extinction time, with the convention $\inf\emptyset=\infty$. A probability measure $\mu$ on $(0,\infty)$ is a QSD if, for all $t\ge0$ and any measurable set $A\subset(0,\infty)$,
\[
\mu(A)=\mathbf P_\mu\bigl(X_t\in A\mid\tau_0^->t\bigr),
\]
where $\mathbf P_\mu$ denotes the distribution of the process $X$ starting from $\mu$.

In applying the QSD criterion below, we shall also use weak continuity of the killed path law with respect to the initial state. This trajectory Feller input is verified later in Theorems \ref{trf} and \ref{trf0} under the same hypotheses.

\btheorem[Quasi-stationary distribution]\label{qsd}
Suppose that $R_0(0)=0$, and that \hyperlink{cond:C4}{\textup{(C4)}} and \hyperlink{cond:C5}{\textup{(C5)}} are satisfied. Assume in addition that either \hyperlink{cond:C1}{\textup{(C1)}} holds, or both \hyperlink{cond:C2}{\textup{(C2)}} and \hyperlink{cond:C3}{\textup{(C3)}} hold. Then there exist $\lambda>0$ and a unique quasi-stationary distribution $\mu$ such that
$$\left\|\mathbf P_\nu\bigl(X_t\in \cdot\mid\tau_0^->t\bigr)-\mu\right\|_{\mathrm{TV}}\le C(\nu)\mrm e^{-\lambda t}, \quad t \ge 0, \ \nu\in \mathcal P((0,\infty)),$$
where the constant $C(\nu) > 0$ depends on  $\nu$.
\etheorem

The proof of Theorem \ref{qsd} is based on the QSD criterion of \cite[Theorem 2.2]{GNW20}. In the verification of this criterion, the Lyapunov input is supplied by Proposition \ref{lyp}, while the trajectory Feller and boundary-continuity inputs are supplied by Theorems \ref{trf}, \ref{trf0}, and \ref{sf0}.

In contrast, when the boundary at zero is non-absorbing ($R_0(0)>0$), the process is prevented from staying at the origin and instead exhibits a long-term persistence within the entire state space $[0, \infty)$. Under the previously established strong Feller and irreducibility properties, such dynamics typically lead to the existence of a unique invariant probability measure. In this non-absorbing regime, we aim to establish the uniform exponential ergodicity, which provides a fast and uniform rate of convergence to the equilibrium regardless of the initial state.
\btheorem[Uniform exponential ergodicity]\label{thmergo}
Suppose that $R_0(0)>0$ and that \hyperlink{cond:C4}{\textup{(C4)}} is satisfied. Assume in addition that either \hyperlink{cond:C1}{\textup{(C1)}} holds, or both \hyperlink{cond:C2}{\textup{(C2)}} and \hyperlink{cond:C3}{\textup{(C3)}} hold. Then there exist constants $C, \lambda>0$ and a unique invariant probability measure $\mu^*$ such that
$$\|\nu P_t-\mu^*\|_{\mathrm{TV}}\le C\mrm{e}^{-\lambda t},\quad t\ge 0,\ \nu\in \mathcal P([0,\infty)).$$
\etheorem

The proof of Theorem \ref{thmergo} uses the Meyn--Tweedie framework in the form of \cite{DMT95}. The additional \(\psi\)-irreducibility, petite-set, and aperiodicity inputs required there are established in Lemma \ref{psiirr}, Proposition \ref{petite}, and Lemma \ref{aperiodic}, respectively.

For reference, the dependence of the main conclusions on the structural
conditions is summarized as follows.
\begin{center}
\small
\begin{tabular}{p{0.41\textwidth}p{0.51\textwidth}}
\hline
Conclusion & Conditions\\
\hline
Strong Feller property & \hyperlink{cond:C1}{\textup{(C1)}} or
\hyperlink{cond:C2}{\textup{(C2)}}\\[2pt]
Irreducibility & \hyperlink{cond:C1}{\textup{(C1)}} or
\hyperlink{cond:C3}{\textup{(C3)}}\\[2pt]
QSD exponential convergence &
\(R_0(0)=0\), \hyperlink{cond:C4}{\textup{(C4)}},
\hyperlink{cond:C5}{\textup{(C5)}},\newline
and either
\hyperlink{cond:C1}{\textup{(C1)}} or both
\hyperlink{cond:C2}{\textup{(C2)}} and
\hyperlink{cond:C3}{\textup{(C3)}}\\[2pt]
Uniform exponential ergodicity &
\(R_0(0)>0\), \hyperlink{cond:C4}{\textup{(C4)}},\newline
and either
\hyperlink{cond:C1}{\textup{(C1)}} or both
\hyperlink{cond:C2}{\textup{(C2)}} and
\hyperlink{cond:C3}{\textup{(C3)}}\\
\hline
\end{tabular}
\end{center}

\begin{exa1} For \(a>0,p>1\) and \(\beta\ge0\), let
\[
R_1(x)=2cx,\qquad R_2(x)=x,\qquad R_0(x)=\beta+b x-a x^p,\qquad x\ge0.
\] The model describes a continuous-state branching process with competition and immigration. The parameter \(\beta\) represents immigration, and \(g(x)=a x^p\) is the competition mechanism. By letting $r=p-1,$ one can find constants \(\lambda_1>0\) and \(\lambda_2\ge0\) such that
\[
R_0(x)=\beta+bx-a x^p
\le -\lambda_1x^p+\lambda_2,\qquad x>0.
\]
Thus Condition \hyperlink{cond:C4}{\textup{(C4)}} holds. When $\beta>0,$ we have $R_0(0)>0.$ Moreover, if \(c>0\), then \hyperlink{cond:C1}{\textup{(C1)}} holds. Alternatively, if there exist constants \(\gamma>0\) and \(\alpha\in(1,2)\) such that
\[
\pi(\d z)\ge
\gamma z^{-1-\alpha}\mathbf 1_{\{0<z\le1\}}\,\d z,
\]
then both \hyperlink{cond:C2}{\textup{(C2)}} and \hyperlink{cond:C3}{\textup{(C3)}} hold. In either case, all the conditions of Theorem \ref{thmergo} are satisfied. Consequently,
the process is uniformly exponentially ergodic, which is  consistent with \cite[Theorem 1.1]{LW20EJP}. 

When $\beta=0,$ the model reduces to a continuous-state branching process with competition. This model is 
introduced by Berestycki et al. \cite{BFF18}. In this case $R_0(0)=0$ and
the boundary point \(0\) is absorbing. Condition \hyperlink{cond:C5}{\textup{(C5)}} is also satisfied: if \(c>0\), this follows from case (i), while in the pure-jump case covered by the above stable-like lower bound on \(\pi\), one may choose \(\theta\in(2-\alpha,1)\), which is non-empty since \(\alpha\in(1,2)\), and then take \(\delta>0\) sufficiently small. Hence the conditions of Theorem \ref{qsd} are satisfied. Thus, there exists a unique quasi-stationary distribution with exponential convergence. This result is consistent with \cite[Theorem 1.2]{lwz}. 
\end{exa1}
\section{Strong Feller and trajectory Feller properties}
\subsection{Localized refined coupling}
We use the refined basic coupling introduced by Luo and Wang \cite{LW18}.
We first recall the coupling generator needed below.

For integrals with respect to the jump measure, and for the corresponding
Poisson integrals, we use the convention
\[
\int_y^x=-\int_x^y=\int_{(y,x]},\qquad
\int_x^\infty=\int_{(x,\infty)},\qquad x\ge y.
\]

Applying It\^{o}'s formula to \eqref{eq:sde_intro}, we deduce that the infinitesimal generator $L$ of the CSNBP is given by
\beqlb\label{gen}
L f(x)
\ar=\ar R_0(x)f'(x)
+\frac12R_1(x)f''(x)\cr
\ar\ar+ R_2(x)\int_0^\infty
\left[f(x+z)-f(x)-zf'(x)\right]\pi(\d z),\qquad f\in C^2_b([0, \infty)).
\eeqlb
In \cite[Proposition 2.2]{LW20EJP}, a two-dimensional Markov coupling process
$\bigl((X_t,Y_t)\bigr)_{t\ge0}$ on $D:=\{(x,y):x\ge y\ge 0\}$ is constructed,
where both $(X_t)_{t\ge0}$ and $(Y_t)_{t\ge0}$ are Markov processes with the same
transition semigroup $(P_t)_{t\ge0}.$ Moreover, $X_t=Y_t$ for all
$t\ge \tilde T,$ where $\tilde T=\inf\{t>0:X_t=Y_t\}$ is the coupling time. We
next give a characterization for the generator of the Markov coupling process
$\bigl((X_t,Y_t)\bigr)_{t\ge0}.$ Let $\Delta=\{(z,z):z\ge 0\}\subset D$ and
$\Delta^c=D\setminus\Delta.$ Given a function $F$ on $D$ twice continuously
differentiable on $\Delta^c,$ we write for any $(x,y)\in\Delta^c,$
{\small\beqlb\label{tildel}
\tilde{L}F(x,y)\ar:=\ar R_0(x)F'_x(x,y)+\frac{1}{2}R_1(x)F''_{xx}(x,y)+R_0(y)F'_y(x,y)+\frac{1}{2}R_1(y)F''_{yy}(x,y)\cr\cr
\ar+\ar  \frac{1}{2} R_{2}(y)\int_{0}^{\infty}\left[F(x+z,x+z)-F(x+z,y+z)\right]\pi_{y-x}(\d z)-\sqrt{R_1(x)R_1(y)}F''_{xy}(x,y)\cr
\ar+\ar \frac{1}{2} R_{2}\left(y\right)\int_{0}^{\infty}\left[F(x+z,2y+z-x)-F(x+z,y+z)\right] \pi_{x-y}(\d z)  \cr\cr
\ar+\ar R_{2}\left(y\right)\int_{0}^{\infty}\left[F(x+z,y+z)-F(x,y)-zF'_x(x,y)-zF'_y(x,y)\right] \pi(\d z) \cr\cr
\ar+\ar ({R_{2}(x)}-{R_{2}(y)})\int_{0}^{\infty} \left[F(x+z,y)-F(x,y)-zF'_x(x,y)\right] \pi(\d z),
\eeqlb}
where $
\pi_{x}=\left(\pi \wedge\left(\delta_{x} * \pi\right)\right),$ for any $x\in\mathbb R.$ Here \(\pi_x\) denotes the common part of \(\pi\) and its shift by \(x\). In \eqref{tildel}, the first line gives the drift and diffusion contributions of the two marginals; the next two lines are the refined basic coupling terms for the common jump parts; and the last two lines correspond to synchronous jumps and to the excess jump intensity \(R_2(x)-R_2(y)\), respectively.

Let $\mathcal D(\tilde L)$ denote the linear space consisting of the functions $F$ such that the integrals in \eqref{tildel} are finite and $\tilde L F$ is locally bounded 
on compact subsets of $\Delta^c.$ We call $(\tilde L,\mathcal D(\tilde L))$ the {\it coupling generator} of the process $(X_t)_{t\ge0}.$ Let 
\[\zeta_n^*=\inf\{t\ge 0:X_t\ge n  \text{\,\,or\,\,} X_t-Y_t\le 1/n\}\]
and assume that $(X_0,Y_0)=(x,y)\in\Delta^c.$ Then, it follows from \cite[Proposition 2.2]{LW20EJP} that for any $n\ge 1$ and $F\in\mathcal D(\tilde L),$
\beqlb\label{eq26043001}
F(X_{t\wedge\zeta_n^*},Y_{t\wedge\zeta_n^*})=F(x,y)+\int_0^{t\wedge\zeta_n^*}\tilde LF(X_s,Y_s)\d s+M_n(t),
\eeqlb
where $(M_n(t))_{t\ge0}$ is a martingale. We write
\(\tilde{\mathbf P}_{x,y}\) and \(\tilde{\mathbf E}_{x,y}\) for the probability
and expectation of this coupling process starting from \((x,y)\).


Note that when $F(x,y)=f(x-y)\I_{\{(x,y)\in\Delta^c\}}$ for some function { $f\in C^2_b([0, \infty)),$} we have $F\in\mathcal D(\tilde L).$ In this case, we also write $\tilde Lf(x,y):=\tilde L F(x,y)$ for convenience. Moreover, for any $(x,y)\in \Delta^c$, 
\beqlb\label{couple1}
\tilde{L}f(x,y)
\ar=\ar
(R_0(x)-R_0(y))f'(x-y)+\frac{1}{2}\left(\sqrt{R_1(x)}+\sqrt{R_1(y)}\right)^2f''(x-y)\cr
\ar\ar +(R_2(x)-R_2(y))\int^\infty_0 [f(x+z-y)-f(x-y)-zf'(x-y)]\pi(\d z)\cr
\ar\ar+\frac{1}{2}R_2(y)[f(2x-2y)-2f(x-y)+f(0)]\pi_{x-y}((0,\infty)),
\eeqlb
where we used the fact that $\pi_{x-y}((0,\infty))=\pi_{y-x}((0,\infty)).$

Recall that $(X_t)_{t\ge0}$ is the solution of \eqref{eq:sde_intro}. Let $\mathbf P_x(\cdot)$ be short for $\mathbf P(\cdot |X_0=x)$ for any $x\in(0,\infty).$ We shall use the following local exit estimate to control short-time departures from a prescribed neighbourhood.

\bglemma\label{lemineq2}
Given two intervals \((a',b')\) and \((a,b)\) such that
\(0<a<a'< b'<b<\infty,\) let $(X_t)_{t\ge0}$ be the solution of \eqref{eq:sde_intro} with initial value $x\in(a',b').$ Define
\[
S=\inf\{t>0:X_t\notin(a,b)\}.
\]
Then there exist constants \(C>0\) and \(t_0>0\), depending on \(a,b,a'\) and \(b'\), such
that, for all \(x\in(a',b')\) and \(0<t\le t_0\),
we have
\[
\mathbf P_x(S<t)\le C\sqrt t.
\]
\edlemma
\begin{proof} For $t>0,$ define
\begin{equation*}
M_{t}=\int_0^{t}\sqrt{R_1({X}_s)}\d B_s+\int_0^{t}\int_0^1\int_0^{R_2({X}_{s-})} z \tilde{N}(\d s,\d z,\d u).
\end{equation*}
Then,
\[M_{t}
=X_{t}- x+X^{(1)}_t+X^{(2)}_t,\quad t\ge 0,  \]
where
\[
X^{(1)}_t=-\int_0^{t}
\left({R_0}({X}_s)-R_2(X_s)\int_1^\infty z\pi(\d z)\right)\d s
\]
and
\[X^{(2)}_t=-\int_0^{t}\int_1^\infty\int_0^{R_2(X_{s-})} z {N}(\d s,\d z,\d u).\]
Let $l=(a'-a)\wedge(b-b').$ For any $x\in (a',b'),$
\begin{equation}\label{3}
\begin{split}
{\mathbf{P}}_{x}\left(S<t\right)
&\leq {\mathbf{P}}_{x}\left(|X_{t\wedge S}-x|\geq l\right)  \\
&\leq\mathbf{ P}_{x}\left(|M_{t\wedge S}|>l/2\right) +\mathbf{ P}_{x}\left(|X^{(1)}_{t\wedge S}|>l/2\right)+\mathbf{ P}_{x}\left(X^{(2)}_{t\wedge S}\neq 0\right).
\end{split}
\end{equation}
Let $M_t^*:=M_{t\wedge S}.$ Then $(M_t^*)_{t\ge0}$ is a martingale whose quadratic variation is
\begin{equation*} \langle M^*\rangle_{t}=\int_0^{t\wedge S}R_1({X_s})\d s+\int_0^{t\wedge S}R_2({X_s}) \d s\int_0^1 z^2\pi(\d z),\quad t\ge 0.\end{equation*}
Using Jensen's inequality and Markov's inequality, we have
\beqlb\label{4}
\mathbf{P}_{x}\left(|M_{t\wedge S}| >l/2 \right)\ar\leq\ar \frac{2}{l} \mathbf{ E}_{x}|M^*_{t}| \leq \frac{2}{l} \left(\mathbf{E}_{x}[(M^*_{t})^2]\right)^{1/2}=\frac{2}{l} \left(\mathbf{E}_{x}\langle M^*\rangle_{t}\right)^{1/2}\cr\cr
\ar=\ar\frac{2}{l}\left\{\mathbf{E}_{x}\left[\int_0^{t\wedge S}R_1({X_s})\d s\right]+\mathbf{ E}_{x}\left[\int_0^{t\wedge S}R_2({X_s}) \d s\int_0^1 z^2\pi(\d z)\right]\right\}^{1/2}\cr\cr
\ar\le\ar \frac{2}{l}\left(\sup_{a\le z\le b}R_1(z)+R_2(b)\int_0^1 z^2\pi(\d z)\right)^{1/2}t^{1/2},
\eeqlb where the last inequality follows from the monotonicity of $R_2.$ Since
\(\int_1^\infty z\pi(\d z)<\infty\), we have
\[
|X^{(1)}_{t\wedge S}|\le t\sup_{a\le z\le b}
\left|R_0(z)-R_2(z)\int_1^\infty y\pi(\d y)\right|.
\]
Hence, for small enough $t>0,$
\begin{equation}\label{6}
\mathbf{P}_{x}\left(|X^{(1)}_{t\wedge S}|
>\frac{l}{2}\right)
\le\I_{\left\{t\sup_{a\le z\le b}
\left|R_0(z)-R_2(z)\int_1^\infty y\pi(\d y)\right|>\frac{l}{2}\right\}}=0.
\end{equation}
Since $R_2$ is non-decreasing and $(X^{(2)}_t)_{t\ge0}$ is a pure jump process, whose waiting time before the first jump can be controlled by an exponential random variable, we have 
\begin{equation}\label{5}
\mathbf{ P}_{x}\left(X^{(2)}_{t\wedge S}\neq 0\right)\le
1-\exp\left(-R_2(b) \pi((1,\infty))t\right).
\end{equation}
Finally, substituting  \eqref{4}, \eqref{6} and \eqref{5} into \eqref{3} gives the desired result.
\qed
\end{proof}

We now formulate the localized part of the coupling argument. A global coupling estimate would require a contraction bound on the whole off-diagonal state space \(\Delta^c\), with constants independent of the spatial position. Such an estimate is not available here because the jump coefficient is discontinuous and the noise coefficients degenerate near the boundary. Instead, after fixing an interior point \(x_0>0\), we work inside a small neighbourhood of \(x_0\) and compare the coupling time with the corresponding local exit time. Condition \textup{(C*)} below is the local Lyapunov condition for the coupling generator that implements this idea. The same local estimate will also be used later to prove the trajectory Feller property.

Here and below, \(b\mcr B(E)\) denotes the set of bounded Borel measurable
functions on a measurable space \(E\).

\hypertarget{cond:C*}{\textbf{(C*)}} For any given $x_0>0$, there exist constants $\delta\in (0,x_0/2)$ and $C>0$ depending on $x_0$, together with a nonnegative function \(F\in\mathcal D(\tilde L)\), such that \(F(x,y)>0\) for \((x,y)\in\Delta^c\),
\[
\sup_{\substack{x,y\in(x_0-\delta,x_0+\delta)\\0<x-y\le r}}F(x,y)\to0
\qquad\text{as }r\downarrow0,
\]
and
\beqlb\label{LY1}
\tilde L F(x,y)\le -C,\qquad
x,y\in(x_0-\delta, x_0+\delta),\quad x>y.
\eeqlb

\bgtheorem\label{11}
Suppose that Condition \hyperlink{cond:C*}{\textup{(C*)}} holds for $x_0$. Then, for any \(t>0\),
\[
\lim_{x\to x_0}\|P_t(x,\cdot)-P_t(x_0,\cdot)\|_{\mathrm{TV}}=0.
\]
In particular, for any \(f\in b\mcr B([0,\infty))\), the map \(x\mapsto P_tf(x)\) is continuous at $x_0$.
\edtheorem

\begin{proof}
Let
\[
T_n=\inf\{t\ge0:X_t-Y_t\le 1/n\},
\]
and
\[
S^*=\inf\{t>0:\text{\,either\,\,}X_t\notin(x_0-\delta, x_0+\delta)\ \text{or}\ Y_t\notin(x_0-\delta, x_0+\delta)\}.
\]
Then the coupling time satisfies
$\tilde T=\lim_{n\to\infty}T_n .$
Let \(\bigl((X_t,Y_t)\bigr)_{t\ge0}\) be the coupling process with generator $\tilde L$ starting from
\((x,y)\in(x_0-\delta/2, x_0+\delta/2)^2\) with $x>y$. Take \(n\) large enough such that \(x-y>1/n\). Then, by \eqref{eq26043001}, we obtain
\beqnn
F(X_{t\wedge T_n\wedge S^*},Y_{t\wedge T_n\wedge S^*})=F(x,y)+\int_0^{t\wedge T_n\wedge S^*}\tilde LF(X_s,Y_s)\d s+M_{t\wedge T_n\wedge S^*}.
\eeqnn

Since \(F\in\mathcal D(\tilde L)\), \(\tilde L F\) is bounded on the compact region where \(x,y\in[x_0-\delta,x_0+\delta]\) and \(x-y\ge1/n\). Hence \((M_{t\wedge T_n\wedge S^*})_{t\ge0}\) is a martingale.
Taking expectations and using the fact that $X_t>Y_t$ for $t\le T_n,$ by \eqref{LY1} we see that
\beqnn
0
&\le&
{\mbb E}
\left[
F(X_{t\wedge T_n\wedge S^*},Y_{t\wedge T_n\wedge S^*})
\right]
\\
&=&
F(x,y)
+
{\mbb E}
\left[
\int_0^{t\wedge T_n\wedge S^*}
\tilde L F(X_s,Y_s)\,\d s
\right]
\\
&\le&
F(x,y)
-
C{\mbb E}
\left[t\wedge T_n\wedge S^*\right].
\eeqnn
Letting \(n\to\infty\) and \(t\to\infty\), we obtain
\beqlb\label{eq:T-S-expectation}
{\mbb E}
\left[\tilde T\wedge S^*\right]
\le
\frac{F(x,y)}{C},
\eeqlb
It follows that, for any \(t>0\),
\beqlb\label{eq:T-S-tail}
{\mbb P}
\left(\tilde T\wedge S^*>t\right)
\le
\frac{F(x,y)}{Ct},
\eeqlb
and
\beqlb\label{eq:S-before-T}
{\mbb P}
\left(F(x,y)^{1/2}<S^*<\tilde T\right)
\le
\frac{F(x,y)^{1/2}}{C}.
\eeqlb

We decompose the event \(\{\tilde T>t\}\) according to whether the localized coupling survives up to time \(t\) or exits the local interval before coupling.
Combining \eqref{eq:T-S-tail} and \eqref{eq:S-before-T}, we obtain
\begin{equation}\label{eq:coupling-prob-main}
\begin{split}
{\mbb P}(\tilde T>t)
&=
{\mbb P}\left(\tilde T\wedge S^*>t\right)
+
{\mbb P}\left(S^*\le t<\tilde T\right)
\\
&\le
{\mbb P}\left(\tilde T\wedge S^*>t\right)
+
{\mbb P}\left(S^*<\tilde T\right)
\\
&\le
{\mbb P}\left(\tilde T\wedge S^*>t\right)
+
{\mbb P}
\left(F(x,y)^{1/2}<S^*<\tilde T\right)
\\
&\quad
+
{\mbb P}
\left(S^*\le F(x,y)^{1/2}\right)
\\
&\le
\frac{F(x,y)}{Ct}
+
\frac{F(x,y)^{1/2}}{C}
+
{\mbb P}
\left(S^*\le F(x,y)^{1/2}\right).
\end{split}
\end{equation}
Define \(S^X=\inf\{t>0:X_t\notin(x_0-\delta, x_0+\delta)\}\) and \(S^Y=\inf\{t>0:Y_t\notin(x_0-\delta, x_0+\delta)\}\).
\[
{\mbb P}
\left(S^*\le F(x,y)^{1/2}\right)\le{\mbb P}
\left(S^X\le F(x,y)^{1/2}\right)+{\mbb P}
\left(S^Y\le F(x,y)^{1/2}\right).
\]
By the above estimate and Lemma \ref{lemineq2}, we obtain
\[
\lim_{x-y\to0+}{\mbb P}
\left(S^*\le F(x,y)^{1/2}\right)=0
\qquad\text{uniformly for }x,y\in(x_0-\delta/2, x_0+\delta/2).
\]
Thus, it follows from
\eqref{eq:coupling-prob-main} that
\beqlb\label{26043002}
{\mbb P}(\tilde T>t)\to0
\qquad\text{uniformly for }x,y\in(x_0-\delta/2, x_0+\delta/2)
\quad\text{as }x-y\to0+.
\eeqlb

For any \(A\in\mcr B([0,\infty))\), applying the coupling inequality gives
\[
\begin{aligned}
|P_t(x,A)-P_t(y,A)|
&=
\left|
{\mbb E}\left[
\bigl(\I_A(X_t)-\I_A(Y_t)\bigr)\I_{\{\tilde T>t\}}
\right]
\right|
\\
&\le
{\mbb P}(\tilde T>t).
\end{aligned}
\]
Taking the supremum over \(A\in\mcr B([0,\infty))\), we obtain
\[
\|P_t(x,\cdot)-P_t(y,\cdot)\|_{\mathrm{TV}}
\le
{\mbb P}(\tilde T>t)\to0
\]
uniformly for \(x,y\in(x_0-\delta/2, x_0+\delta/2)\) as \(x-y\to0+\). The other ordering is handled in the same way. Taking \(y=x_0\) gives the asserted total variation continuity at \(x_0\). The strong Feller property follows immediately.

\qed
\end{proof}

\begin{rem1}
The main advantage of the preceding criterion is that it reduces the
strong Feller property on \((0,\infty)\) to estimates on arbitrary bounded
subintervals. Thus, one does not need a global coupling estimate on the whole state space.
This makes the criterion particularly useful for processes with boundary
degeneracy and state-dependent jump intensities.
\end{rem1}

\subsection{Strong Feller property}
\bgproposition\label{sfc1}
Condition \hyperlink{cond:C*}{\textup{(C*)}} is satisfied if either
\hyperlink{cond:C1}{\textup{(C1)}} or
\hyperlink{cond:C2}{\textup{(C2)}} holds.
\edproposition

\begin{proof}
For \(\rho\in(0,1/2]\), define
\beqlb\label{phi}
\phi(r)=1-\mrm e^{-r^\rho},\qquad r\ge0.
\eeqlb
For any $r>0$, \[\phi(0)=0,\quad\phi'(r)=\rho r^{\rho-1} \mrm{e}^{-r^\rho}>0, \quad
\phi''(r)=-\rho r^{\rho-2} \mrm{e}^{-r ^\rho}[1-\rho+\rho r^\rho]<0,\]
and
$$\phi^{(3)}(r)=\rho(2-\rho)r^{\rho-3}\mrm{e}^{-r ^\rho}[1-\rho+\rho r^\rho]+\rho^2r^{2\rho-3}\mrm{e}^{-r ^\rho}[1-2\rho+\rho r^\rho]\ge0.$$

Set \(F(x,y)=\phi(x-y)\I_{\{(x,y)\in\Delta^c\}}\).
Since the possible singularity of \(F\) is confined to the diagonal, \(F\in\mathcal D(\tilde L)\).
Recalling the identity \eqref{couple1}, we have, for $0<y<x,$
\begin{equation}\label{revised2}\begin{split}
 \tilde L F(x,y)= &\left(R_0(x)-R_0(y)\right)\phi'(x-y)+\frac{1}{2}\left(\sqrt{R_1(x)}+\sqrt{R_1(y)}\right)^2\phi''(x-y)\\
 &+\left(R_2(x)-R_2(y)\right) \int_0^\infty
\left(\phi(x-y+z)-\phi(x-y)-\phi'(x-y)z\right)\,\pi(\d z)\\
 &+ \frac{1}{2} R_2(y)\left[\phi(2(x-y)
)-2\phi(x-y)\right]\pi_{x-y}((0,\infty))\\
=:&I_1+I_2+I_3+I_4.
\end{split}\end{equation}
By the mean value theorem and the facts that $\phi(0)=0$, $\phi''\le 0$ and $\phi^{(3)}\ge0$,
we obtain for $x>0,$
\beqlb\label{es4281}\phi(2x)-2\phi(x)\ar=\ar\int_0^{2x}\phi'(u)\,\d u-2\int_0^x\phi'(u)\,\d u=\int_0^x\phi'(x+u)\,\d u-\int_0^x\phi'(u)\,\d u\cr\cr
	\ar\le\ar\int_0^x \phi''(2x)x\,\d u=\phi''(2x)x^2<0.\eeqlb
Hence
\(I_4\le0.\) And \(I_2\le0\) since \(\phi''\) is negative on \((0,\infty)\). In addition, since \(R_2\) is non-decreasing, we also
have \(I_3\le0\). Furthermore, by the assumption \eqref{lclip}, there is a constant
\(k_0>0\) such that
\[
|R_0(x)-R_0(y)|\le k_0(x-y),\qquad x,y\in(x_0-\delta, x_0+\delta),\ x>y.
\]
Hence
\[
I_1\le k_0\rho(x-y)^\rho \mrm e^{-(x-y)^\rho}\le 2k_0\rho\delta^\rho.
\]
Thus, $I_1$ is bounded on $(x_0-\delta, x_0+\delta).$ Therefore, \hyperlink{cond:C*}{\textup{(C*)}} holds if one of $I_2,I_3,I_4$ tends to $-\infty$ as \(x-y\to0+\).

In the case of \hyperlink{cond:C1}{\textup{(C1)}}, since \(\inf_{z\in(x_0-\delta,x_0+\delta)}R_1(z)>0\) and  \(\rho\in(0,1/2]\),
we have
\[
I_2
=
-\frac{1}{2}\left(\sqrt{R_1(x)}+\sqrt{R_1(y)}\right)^2
\rho (x-y)^{\rho-2}\mrm e^{-(x-y)^\rho}
\left[1-\rho+\rho(x-y)^\rho\right]
\to-\infty
\]
as \(x-y\to0+\). 

In the case of \hyperlink{cond:C2}{\textup{(C2)}},
\cite[Example 1.2]{LW18} gives the following lower bound for the common part of the shifted measures. Indeed, by the lower bound in \hyperlink{cond:C2}{\textup{(C2)}}, if \(0<z<1/4\), then on the interval \((2z,1/2]\), both \(\pi(\d u)\) and \((\delta_z*\pi)(\d u)\) dominate \(\gamma u^{-1-\alpha}\d u\). Hence, after decreasing the constant if necessary, there is \(\hat C>0\) such that
\[
\pi_z(\mathbb R_+)\ge \gamma\int_{2z}^{1/2}u^{-1-\alpha}\d u
\ge \hat Cz^{-\alpha}
\]
for all sufficiently small \(z\), where \(\hat C\) depends on \(\gamma\) and \(\alpha\).
Combining this estimate with \eqref{es4281}, we obtain,
\beqnn
I_4
\ar\le\ar
-\frac{1}{2} \hat CR_2(y)\rho (2x-2y)^{\rho-2}
\mrm e^{-(2x-2y)^\rho}
\left[1-\rho+\rho (2x-2y)^\rho\right](x-y)^{2-\alpha}
\cr\cr\ar\le\ar
-2^{\rho-3} \hat CR_2(y)\rho (x-y)^{\rho-\alpha}
\mrm e^{-(2x-2y)^\rho}
\left[1-\rho+\rho (2x-2y)^\rho\right].
\eeqnn
Since \(R_2\) is strictly positive on \((0,\infty)\), we have $\inf_{y\in(x_0-\delta,x_0+\delta)}R_2(y)>0.$
By taking
\(
\rho<\frac{1}{2}(\alpha\wedge1),
\)
we observe that for $x_0-\delta<y<x<x_0+\delta,$
\[
I_4\to-\infty
\qquad\text{as }x-y\to0+.
\]
In both cases, by decreasing \(\delta\) if necessary, the above convergence to
\(-\infty\) yields the required uniform bound \(\tilde L F(x,y)\le -C\) on the
local region in \hyperlink{cond:C*}{\textup{(C*)}}.
Hence the assertion follows.
\qed
\end{proof}

{\bf Proof of Theorem \ref{SF_a}.}
Using Theorem \ref{11} and Proposition \ref{sfc1}, we have proved that
under \hyperlink{cond:C1}{\textup{(C1)}} or \hyperlink{cond:C2}{\textup{(C2)}}, for each $t>0$ and every bounded Borel function $f$ on $[0,\infty),$ $x\mapsto P_tf(x)$ is continuous on $(0,\infty).$
\qed

\btheorem\label{sf0}
Suppose that \hyperlink{cond:C5}{\textup{(C5)}} is satisfied. Then for any $t>0$ and any
$f\in b\mcr{B}([0, \infty))$,
\[
\lim_{x\to0+}P_tf(x)=P_tf(0).
\]
\etheorem

\begin{proof}
{\bf Step 1.} Let \(\phi\) be defined by \eqref{phi}.  We first show that under \hyperlink{cond:C5}{\textup{(C5)}}, there exists \(\delta\in(0,1)\) such that
\[
L\phi(x)\le -1,\qquad 0<x<\delta.
\]
It is sufficient to show that
\[
\lim_{x\to0+}L\phi(x)=-\infty.
\]
For \(0<x<\delta\), by the definition of \(L\) in \eqref{gen}, we have
\beqnn
L\phi(x)
\ar=\ar
R_0(x)\phi'(x)
+\frac12R_1(x)\phi''(x)
+
R_2(x)\int_0^\infty
\left[\phi(x+z)-\phi(x)-z\phi'(x)\right]\pi(\d z)
\cr
\ar=:\ar J_1+J_2+J_3.
\eeqnn
By \hyperlink{cond:C5}{\textup{(C5)}},
\[
J_1\le \rho kx^{\theta+\rho-1}\mrm e^{-x^\rho}:=J(x).
\]
Since \(\phi\) is concave, we have \(J_2,J_3\le0\). 

In the case of \hyperlink{cond:C5}{\textup{(C5)}}\textup{(i)}, choosing
\(\rho<1-\theta\), we have
$$J_2\le -\frac{2k}{\theta}\frac12\rho(1-\rho)\mrm e^{-x^\rho}x^{\theta+\rho-1}=-\frac{k}{\theta}\rho(1-\rho)\mrm e^{-x^\rho}x^{\theta+\rho-1}= -c_2J(x),$$
where $c_2=\frac{1-\rho}{\theta}>1$.

In the case of \hyperlink{cond:C5}{\textup{(C5)}}\textup{(ii)}, choose \(0<\rho<\frac12\wedge(1-\theta).\) 
Since the integrand in \(J_3\) is non-positive, for $x\in(0,1),$ we have
\beqnn
J_3
\ar\le\ar
R_2(x)\int_0^x
\left[\phi(x+z)-\phi(x)-z\phi'(x)\right]\pi(\d z)\le\frac12R_2(x)\phi''(2x)\int_0^xz^2\pi(\d z)
\cr\cr\ar\le\ar
\frac{6(2-\alpha)}{\gamma}kx^{\theta+\alpha-1}\phi''(2x)\frac{\gamma}{2-\alpha}x^{2-\alpha}\cr\cr
\ar=\ar-\rho kx^{\theta+\rho-1}6\cdot2^{\rho-2}\mrm e^{-(2x)^\rho}
\left[1-\rho+\rho(2x)^\rho\right].
\eeqnn
Hence, for small enough $x>0$,
\beqnn
J_3\le -c_1\rho kx^{\theta+\rho-1}\mrm e^{-(2x)^\rho}= -c_1J(x)(\e^{x^\rho-(2x)^\rho}),
\eeqnn
where \(c_1:=6\cdot2^{\rho-2}(1-\rho)>1\) after decreasing \(\rho\) if necessary.

In either case, after reducing \(\delta\) if necessary, there exists a constant
\(c_*>1\) such that
\[
J_2+J_3\le -c_*J(x),\qquad 0<x<\delta.
\]
Therefore, under \hyperlink{cond:C5}{\textup{(C5)}}, for small $x$,
$$L\phi(x)=J_1+J_2+J_3\le -(c_*-1)J(x)\to-\infty,\text{\,\,\,as\,\,\,}x\to 0+.$$

{\bf Step 2.} Since \(R_0(0)=R_1(0)=R_2(0)=0\), the state \(0\) is absorbing.
Now, let \(x\in(0,\delta)\). For \(n\ge1\), let
\[
\tau_\epsilon^-=\inf\{t\ge0:X_t\le \epsilon\},\quad \tau_\delta^+=\inf\{t\ge0:X_t\ge \delta\},\quad\text{for\,\,} \delta>\epsilon\ge 0.
\]
Similar to the proof of Theorem \ref{11}, by It\^o's formula we have
\[\phi(X_{t\wedge\tau_\epsilon^-\wedge\tau_\delta^+})=\phi(x)+\int_0^{t\wedge\tau_\epsilon^-\wedge\tau_\delta^+}L\phi(X_s)\d s+M_t,\]
where $(M_t)_{t\ge0}$ is a martingale. Taking expectations on both sides of the above, we obtain
\beqlb\label{es04291}
0
\le
\mathbf E_x[\phi(X_{t\wedge\tau_\epsilon^-\wedge\tau_\delta^+})]
\le
\phi(x)-
\mathbf E_x[t\wedge\tau_\epsilon^-\wedge\tau_\delta^+].
\eeqlb
Letting \(t\to\infty\) in the inequality above, we obtain
\[
\mathbf E_x[\tau_\epsilon^-\wedge\tau_\delta^+]\le \phi(x).
\]
Since $\phi$ is continuous and $\phi(0)=0,$ using Markov's inequality, for any $t>0,$ we have
\[\mathbf P_x\left(\tau_0^-\wedge\tau_\delta^+>t\right)\le\phi(x)/t\rightarrow 0\text{\,\,\,as\,\,\,}x\rightarrow 0.\]
On the other hand, since $\phi$ is increasing, using \eqref{es04291} again gives
\[\mathbf P_x\left(\tau_\delta^+\le t\wedge\tau_0^-\right)\le\phi(x)/\phi(\delta)\rightarrow 0\text{\,\,\,as\,\,\,}x\rightarrow 0.\]

Combining the two inequalities above gives
\beqlb\label{es04292}\mathbf P_x\left(\tau_0^->t\right)\le\mathbf P_x\left(\tau_0^-\wedge\tau_\delta^+>t\right)+\mathbf P_x\left(\tau_\delta^+\le t\wedge\tau_0^-\right)\le\phi(x)/t+\phi(x)/\phi(\delta),\eeqlb
which tends to $0$ as $x\to 0.$ Thus, for any $f\in b\mcr B([0,\infty))$ and $t>0,$
\beqnn
|P_tf(x)-P_tf(0)|\ar\le\ar \left|\mathbf E_x\left[\left(f(X_t)-f(0)\right)\I_{\{\tau_0^->t\}}\right]\right|\cr\cr
\ar\le\ar 2\|f\|_\infty\mathbf P_x(\tau_0^->t)\to 0 \text{\,\,\,as\,\,\,} x\to 0.
\eeqnn
The proof is complete.
\qed
\end{proof}

\subsection{Trajectory Feller property}

We now establish the trajectory Feller input needed for the QSD application. The criterion in \cite{GNW20} requires weak continuity of the killed path law with respect to the initial condition.

For
$t_0>0,$ denote by $D([0,t_0],\mathbb R_+)$ the space of $\mathbb R_+$-valued c\`adl\`ag paths defined on $[0,t_0]$ and equipped with the Skorokhod topology.

\btheorem\label{trf}
Suppose that \hyperlink{cond:C1}{\textup{(C1)}} or \hyperlink{cond:C2}{\textup{(C2)}} holds. Then for any $t_0>0,$ $$z\mapsto \mathbf P_z\left((X_t)_{t\in[0,t_0]}\in\cdot\right)$$ is continuous on $(0,\infty)$ in the sense of weak convergence.\etheorem

\begin{proof}
Fix $t_0>0$ and $x_0\in (0,\infty)$. We prove that $z\mapsto \mathbf P_z\left((X_t)_{t\in[0,t_0]}\in\cdot\right)$ is continuous at $z=x_0$. Let $\bigl((X_t,Y_t)\bigr)_{t\ge0}$ be the coupling Markov process with generator $\tilde L$ with $X_0=x$ and $Y_0=y.$ It suffices to consider the ordered case \(y<x\), since the case \(x<y\) follows by relabeling the two initial points. Noting that the Skorokhod distance is dominated by the supremum norm, it is sufficient to prove that for any fixed $t_0,x_0, \varepsilon_1,\varepsilon_2>0$, there exists a small enough $\delta>0$ such that for any  $0<x_0-\delta<y<x<x_0+\delta,$
\beqnn
\tilde{\mathbf{P}}_{x,y}\left(\sup_{0 \le t \le t_0} (X_t - Y_t)>\varepsilon_1\right)<\varepsilon_2.
\eeqnn 
Let
\[T_X= \inf\{t
\ge 0:|X_t-x_0|>\varepsilon_1/6\}, T_Y=\inf\{t\ge 0:|Y_t-x_0|>\varepsilon_1/6\}.\]
According to Lemma \ref{lemineq2}, there exist $t_1\in(0,t_0)$ and 
$\delta\in\big(0,\varepsilon_1/6\wedge x_0\big)$ such that for any 
$x_0-\delta<y<x<x_0+\delta$,
\[\tilde{\mathbf P}_{x,y}\left(T_X<t_1\right)<\varepsilon_2/3\quad\text{and}\quad \tilde{\mathbf P}_{x,y}\left(T_Y<t_1\right)<\varepsilon_2/3.\]
This yields the result that for any $0<x_0-\delta<y<x<x_0+\delta$ with $x-y<\varepsilon_1/3,$
\beqlb\label{10311}
\ar\ar\tilde{\mathbf{P}}_{x,y}\left(\sup_{0 \le t \le t_1} (X_t - Y_t)>\varepsilon_1\right)\cr
\ar\le\ar \tilde{\mathbf{P}}_{x,y}\left(\sup_{0 \le t \le t_1} |X_t - x|>\varepsilon_1/3\right)+\tilde{\mathbf{P}}_{x,y}\left(\sup_{0 \le t \le t_1} |Y_t-y|>\varepsilon_1/3\right)+\I_{\{(x - y)>\varepsilon_1/3\}}\cr
\ar\le\ar
\tilde{\mathbf P}_{x,y}\left(T_X<t_1\right)+\tilde{\mathbf P}_{x,y}\left(T_Y<t_1\right)\le 2\varepsilon_2/3.
\eeqlb
On the other hand, by Proposition \ref{sfc1}, Condition \hyperlink{cond:C*}{\textup{(C*)}} is satisfied if either
\hyperlink{cond:C1}{\textup{(C1)}} or
\hyperlink{cond:C2}{\textup{(C2)}} holds. Then, applying \eqref{26043002}, there exists small enough $\delta\in (0,\varepsilon_1/6)$ such that for any $0<x_0-\delta<y<x<x_0+\delta,$
\beqlb\label{eq22d}
\tilde{\mathbf{P}}_{x,y}(\tilde T>t_1)<\varepsilon_2/3,
\eeqlb
where $\tilde T=\inf \left\{t \geq 0: X_{t}=Y_{t}\right\}$ is the coupling time.
Combining \eqref{10311} and \eqref{eq22d}, we have, for any fixed $t_0>0$ and any fixed $x_0\in (0,\infty)$, if $0<x_0-\delta<y<x<x_0+\delta,$ then
\beqlb\label{10161}
\ar\ar\tilde{\mathbf{P}}_{x,y}\left(\sup_{0 \le t \le t_0} (X_t - Y_t)>\varepsilon_1\right)\cr
\ar=\ar\tilde{\mathbf{P}}_{x,y}\left(\sup_{0 \le t \le t_0} (X_t - Y_t)>\varepsilon_1, \tilde T> t_1\right)+\tilde{\mathbf{P}}_{x,y}\left(\sup_{0 \le t \le t_1} (X_t - Y_t)>\varepsilon_1, \tilde T\le t_1\right)\cr
\ar\le\ar
\tilde{\mathbf{P}}_{x,y}\left(\tilde T> t_1\right)+\tilde{\mathbf{P}}_{x,y}\left(\sup_{0 \le t \le t_1} (X_t - Y_t)>\varepsilon_1\right)<\varepsilon_2.
\eeqlb
 This proves that the mapping $z \mapsto \mathbf{P}_z\left((X_t)_{t \in [0,t_0]} \in \cdot\right)$ is continuous on $(0, \infty)$ in the sense of weak convergence in the Skorokhod topology. 
\qed
\end{proof}

\bgproposition[Compact-starting finite-time boundedness]\label{prop_compact_bounded}
Suppose that \hyperlink{cond:C1}{\textup{(C1)}} or \hyperlink{cond:C2}{\textup{(C2)}} holds. Let \(K\subset(0,\infty)\) be compact and let \(T>0\). If
\[
\tau_R^+=\inf\{t\ge0:X_t\ge R\},
\]
then
\[
\lim_{R\to\infty}\sup_{x\in K}\mathbf P_x(\tau_R^+\le T)=0.
\]
\edproposition
\proof
By Theorem \ref{trf}, the map
\[
x\mapsto \mathbf P_x\bigl((X_t)_{t\in[0,T]}\in\cdot\bigr)
\]
is continuous on \((0,\infty)\) in the sense of weak convergence. Hence its image over the compact set \(K\) is compact in the weak topology of probability measures on \(D([0,T],\mathbb R_+)\). Since \(D([0,T],\mathbb R_+)\) is a Polish space, Prokhorov's theorem implies that this compact family is tight. Thus, for any \(\varepsilon>0\), choose a compact set \(\Gamma_\varepsilon\subset D([0,T],\mathbb R_+)\) such that
\[
\inf_{x\in K}\mathbf P_x\bigl((X_t)_{t\in[0,T]}\in\Gamma_\varepsilon\bigr)>1-\varepsilon.
\]
Every compact subset of \(D([0,T],\mathbb R_+)\) is bounded in the supremum norm. Indeed, the map \(\omega\mapsto\sup_{0\le t\le T}\omega(t)\) is continuous under the Skorokhod topology. Hence there exists \(R_\varepsilon>0\) such that
\[
\sup_{0\le t\le T}\omega(t)<R_\varepsilon,\qquad \omega\in\Gamma_\varepsilon.
\]
For \(R\ge R_\varepsilon\), this gives
\[
\sup_{x\in K}\mathbf P_x(\tau_R^+\le T)\le\varepsilon.
\]
The result follows.
\qed

\bgtheorem\label{trf0}
Suppose that \hyperlink{cond:C5}{\textup{(C5)}} holds. Then for any $t_0>0$, 
\[
z\mapsto \mathbf P_z\bigl((X_t)_{t\in[0,t_0]}\in\cdot\bigr)
\]
is right-continuous at $z=0$ in the sense of weak convergence.
\edtheorem

\begin{proof}
The proof is similar to that of Theorem \ref{trf}. Since the Skorokhod distance
is dominated by the supremum norm, it is enough to show that for any fixed
$t_0,\varepsilon_1,\varepsilon_2>0$, there exists $\delta>0$ such that for all
$x\in(0,\delta)$,
\[
\mathbf P_x\left(\sup_{0\le t\le t_0}X_t>\varepsilon_1\right)<\varepsilon_2.
\]
 Let $\tau^+_{\varepsilon_1}=\inf\{t\ge0:X_t\ge\varepsilon_1\}.$
Applying the proof of Lemma \ref{lemineq2} to the one-sided interval
\([0,\varepsilon_1]\), and using the boundedness of the coefficients on this
interval together with \(\varepsilon_1-x\ge \varepsilon_1/2\), there exists a constant
$C>0$ such that, for $ 0\le x\le \varepsilon_1/2$ and small enough $t>0,$
\[
\mathbf P_x(\tau^+_{\varepsilon_1}<t)\le C\sqrt t.\]
Hence, for  small enough \(t_1>0\), we have
\[
\mathbf P_x\left(\sup_{0\le t\le t_1}X_t>\varepsilon_1\right)<\varepsilon_2/2,\qquad \forall x\in(0, \varepsilon_1/2)
\]
By Theorem \ref{sf0}, applied to \(f=\I_{(0,\infty)}\), and since \(0\) is absorbing, taking $\delta\in(0,\varepsilon_1)$ small enough, we have
\[
\mathbf P_x(\tau_0^- > t_1) < \varepsilon_2/2, \qquad x\in(0,\delta).
\]
Since $X_t=0$ for all $t\ge\tau_0^-$, combining the above two estimates, for any $x\in(0,\delta)$,
\[
\mathbf P_x\Bigl(\sup_{0\le t\le t_0}X_t > \varepsilon_1\Bigr)
\le \mathbf P_x(\tau_0^- > t_1) + \mathbf P_x\Bigl(\sup_{0\le t\le t_1}X_t > \varepsilon_1\Bigr)
< \varepsilon_2.
\]
This proves the right-continuity at $0$ in the Skorokhod topology.
\qed
\end{proof}

\section{Irreducibility}
Establishing irreducibility for the non-linear process \eqref{eq:sde_intro} is highly non-trivial. Standard tools, such as the Girsanov measure transformation, heavily rely on the uniform ellipticity of the diffusion term. In our setting, the diffusion coefficient $R_1(x)$ degenerates at zero, rendering these classical methods inapplicable. Furthermore, unlike linear continuous-state branching processes where a global Lamperti-type time-change can fully decouple the system into a space-homogeneous L\'evy process, our non-linear state-dependent rates ($R_0, R_1, R_2$) prevent such simultaneous diagonalization.

To overcome these intrinsic difficulties, we adopt a purely probabilistic approach based on state-dependent time-change, jump truncation, and pathwise coincidence. The proof is organized around the following implication:
\[
\begin{aligned}
\text{upward reachability}+\text{downward reachability}
&\Longrightarrow \text{hitting every interior level}  \\
&\Longrightarrow \text{irreducibility on open sets}.
\end{aligned}
\]
\begin{itemize}
    \item \textbf{Upward reachability:} We construct an auxiliary process that truncates large jumps. By isolating the jump mechanism and exploiting the spatial independence of Poisson random measures, we show that an upward jump event with positive probability occurs, during which the true process coincides with the auxiliary one.
    \item \textbf{Downward reachability:} Driving the process downward is particularly challenging since the branching mechanism admits only positive jumps. We first apply a state-dependent time-change driven by $R_2$ to normalize the jump intensity. We then introduce an auxiliary process $Y^\delta$ in which jumps larger than $\delta$ are suppressed while their negative compensator is retained. Under the infinite-variation condition ($\int_0^1 z\pi(\d z)=\infty$), the small-jump part of this compensator gives a negative drift, which becomes arbitrarily large as $\delta\downarrow0$ and pushes the path toward lower levels. In the finite-variation case, the same downward passage is instead supplied by the non-degenerate interior diffusion. Finally, we show that, with positive probability, the original process coincides with the corresponding auxiliary path.
\end{itemize}

In this section, we focus on the irreducibility of the process. Recall that $(X_t)_{t\ge0}$ is the solution of \eqref{eq:sde_intro}, with its transition semigroup denoted by $(P_t)_{t\ge0}$. Although the compensated form contains a negative compensator, all actual jumps of \(X\) are nonnegative; in particular, downward level crossings occur continuously.

For any process $(\xi_t)_{t\ge0}$, and any constants $0\le\epsilon<x<z<\infty$, we denote the first passage times by 
$$\tau_\epsilon^-(\xi):=\inf\{t\ge 0: \xi_t\le\epsilon\},\quad \tau_z^+(\xi):=\inf\{t\ge 0: \xi_t\ge z\},$$
$$\tau_z^{\epsilon,+}(\xi):=\inf\{t\ge 0: \xi_t \ge z, t < \tau^-_{\epsilon}(\xi)\},\quad \tau_\epsilon^{z,-}(\xi):=\inf\{t\ge 0: \xi_t \le \epsilon, t < \tau^+_{z}(\xi)\},$$
and $\tau_y(\xi):=\inf\{t\ge 0: \xi_t=y\}$ for all $y\ge 0$. Recall that $\mathbf P_x(\cdot)$ is short for $\mathbf P(\cdot \mid X_0=x)$ for any $x\in(0,\infty)$.

\bgproposition\label{prop_reach_irred}
Suppose that for any $t>0$ and any $0<\epsilon<x<y$, $\mathbf P_x(\tau^-_\epsilon(X)<t)>0$ and $\mathbf P_x(\tau^+_y(X)<t)>0$. Then the transition semigroup $(P_t)_{t\ge0}$ is irreducible in $(0, \infty)$; that is, $P_t(x, O) > 0$ for any $x>0, t>0$ and nonempty open set $O \subset (0, \infty)$.
\edproposition
\proof
{\bf Step 1.} We first show that for any $z>0$ and $x>0$, $\mathbf P_x(\tau_z(X)<t)>0$.
First, consider the case $z \le x$. Take \(0<\epsilon<z\). Since the drift and Brownian parts are continuous and the jump part has only positive jumps, the process cannot jump downward across \(z\). Hence, on the event $\{\tau^-_\epsilon(X)<t\}$, the path must hit \(z\) before reaching \(\epsilon\). Therefore,
\[
\mathbf P_x(\tau_z(X)<t) \ge \mathbf P_x(\tau^-_\epsilon(X)<t) > 0.
\]

We next assume $z>x$. Applying the strong Markov property at $\tau_z^+(X)$, we have
\beqnn
\mathbf P_x\left(\tau_z(X)<t\right)\ar\ge\ar\mathbf P_x\left(\tau^+_z(X)<t/2, \, \tau^-_\epsilon(X)\circ\theta_{\tau^+_z(X)}<t/2\right)\cr\cr
\ar=\ar\mathbf E_x\left[\I_{\{\tau^+_z(X)<t/2\}}\mathbf P_{X_{\tau^+_z(X)}}\left(\tau^-_\epsilon(X)<t/2\right)\right],
\eeqnn
where $\theta$ denotes the shift operator, and we choose an arbitrary $\epsilon < z$. According to our hypothesis, the mapping $y \mapsto \mathbf P_y(\tau^-_\epsilon(X)<t/2)$ is strictly positive for all $y \ge z$ and $\mathbf P_x(\tau^+_z(X)<t/2)>0$. Since the process has no negative jumps, any path starting from a point \(y\ge z\) and reaching \(\epsilon<z\) must hit \(z\) along the way. Hence, the above expectation is strictly positive.

{\bf Step 2.} We now prove irreducibility on open sets.
Fix a nonempty open set \(O\subset(0,\infty)\), and pick \(z\in O\). 
By Lemma \ref{lemineq2}, there exists \(h_0>0\) such that
\beqlb\label{eq:local-positivity-at-z}
P_u(z,O)>0,\qquad 0<u\le h_0.
\eeqlb
We now transfer this local positivity at \(z\) to arbitrary initial points. For any \(y>0\) and \(0<u\le h_0\),  applying the strong Markov property at \(\tau_z(X)\), we get
\[
P_u(y,O)
\ge
\mathbf E_y\left[
\I_{\{\tau_z(X)<u\}}
P_{u-\tau_z(X)}(z,O)
\right],
\]
In the event \(\{\tau_z(X)<u\}\), we have \(u-\tau_z(X)\in(0,h_0]\), and hence
\(P_{u-\tau_z(X)}(z,O)>0\) by \eqref{eq:local-positivity-at-z}. Since \(\mathbf P_y(\tau_z(X)<u)>0\) by Step 1, this gives
\beqlb\label{eq:small-time-into-O}
P_u(y,O)>0,\qquad y>0,\quad 0<u\le h_0.
\eeqlb

Now fix any \(x>0\) and \(t>0\). Choose an integer \(n\ge1\) such that \(h:=t/n\le h_0\). To complete the proof, it suffices to show that \(P_{kh}(x,O)>0\) for \(k=1,\ldots,n\). The case \(k=1\) follows directly from \eqref{eq:small-time-into-O}. Then the Chapman-Kolmogorov equation gives
\[
P_{2h}(x,O)
\ge
\int_{(0,\infty)}P_{h}(x,\d y)P_h(y,O)>0.
\]
Thus, by induction, \(P_t(x,O)>0\), proving irreducibility in \((0,\infty)\).
\qed

\subsection{Upward reachability}

We first prove upward reachability for the original process \(X\). This part does not use the time-change; upward movement is supplied either by positive jumps or by the non-degenerate diffusion in the interior.

\bglemma\label{upup}
Suppose that $\pi((0, \infty))>0$. Then for any $y>x>0$, we have $\mathbf{P}_x(\tau^+_y(X)< t)>0$.
\edlemma
\proof
Since $\pi((0, \infty))>0$, there exists some $a>0$ such that $\pi((a,\infty))>0$. For any given target level $y>x$, we can choose an integer $N \ge 1$ such that $x + N(a/2) \ge y$. Thus, it is sufficient to show that for any starting point $x_0>0$ and any given time duration $t_0>0$, the process can move upward by at least $a/2$ with strictly positive probability, i.e. $\mathbf{P}_{x_0}(\tau^+_{x_0+a/2}(X) < t_0) > 0$. By taking $t_0 = t/N$ and applying the strong Markov property successively at the hitting times in the $N$ steps, the general conclusion $\mathbf{P}_x(\tau^+_y(X)< t)>0$ follows immediately. We next focus on proving this single-step movement.

Fix any $x_0>0$ and $t_0>0$, and set \(s_0=t_0/2\) and
\(\ell=(x_0\wedge a)/2\). We decompose the process by isolating jumps of size strictly greater than $a$:
\beqnn
X_{t} \ar=\ar x_0 + \int_0^{t} \widehat{R}_0(X_s) \d s + \int_0^{t} \sqrt{R_1(X_s)} \d B_s + \int_0^{t} \int_0^a \int_0^{R_2(X_{s-})} z \tilde{N}(\d s, \d z, \d u) \cr
\ar\ar + \int_0^{t} \int_a^\infty \int_{0}^{R_2(X_{s-})} z N(\d s, \d z, \d u), 
\eeqnn
where $\widehat{R}_0(\xi) := R_0(\xi) - R_2(\xi) \int_{a}^\infty z \pi(\d z)$.
We introduce an auxiliary process $(Y_t)_{t\ge0}$ starting at $x_0$, defined by removing the large jumps:
\beqnn
Y_{t}\ar=\ar x_0+\int_0^{t} \widehat{R}_0(Y_s)\d s+\int_0^{t}\sqrt{R_1({Y}_s)}\d B_s+\int_0^{t}\int_0^a\int_0^{R_2({Y}_{s-})} z{\tilde{N}}(\d s, \d z, \d u).
\eeqnn
Let $\sigma_1$ denote the arrival time of the first large jump along the auxiliary path:
$$ \sigma_1= \inf\left\{t > 0 : \int_0^{t} \int_a^\infty \int_{0}^{R_2(Y_{s-})} N(\d s, \d z, \d u) > 0\right\}. $$ 
Therefore, by pathwise uniqueness, we have $X_t=Y_t$ for all $t<\sigma_1$. 

Fix \(m>x_0\) and set
\[
S_m=s_0 \wedge \tau^-_{x_0-\ell}(Y)\wedge\tau_m^+(Y).
\]
By the spatial independence of Poisson random measures, conditioned on the path of the auxiliary process, the probability of observing no large jumps before \(S_m\) is
\beqnn
\mathbf{P}_{x_0}\left(\sigma_1 > S_m \mid Y \right)
= \exp\left\{- \pi((a, \infty))\int_0^{S_m} R_2(Y_{s-}) \d s\right\}.
\eeqnn
Since \(x_0-\ell>0\), \(Y_0=x_0\), the paths are right-continuous, and \(R_2\) is strictly positive on \((0,\infty)\), the integral in the exponent is strictly positive almost surely. It is also finite, because \(Y\) is stopped in the compact interval \([x_0-\ell,m]\). Hence, the conditional probability above is strictly smaller than one almost surely. Taking expectations yields $\mathbf{P}_{x_0}(\sigma_1 \le S_m) > 0$. 

On the event $\{\sigma_1 \le S_m\}$, we know that right before the jump, $Y_{\sigma_1-} \ge x_0 - \ell$. Since the isolated jump is of size at least $a$ and \(\ell\le a/2\), we immediately obtain $X_{\sigma_1} = Y_{\sigma_1-} + \Delta X_{\sigma_1} \ge x_0 - \ell + a \ge x_0 + a/2$. Thus, 
$$\mathbf{P}_{x_0}(\tau^+_{x_0+a/2}(X) < t_0) \ge\mathbf{P}_{x_0}(\sigma_1 \le S_m)>0,$$ which completes the proof.
\qed
\bglemma\label{upupri}
Suppose that \hyperlink{cond:C1}{\textup{(C1)}} holds. Then for any $y>x>0$, we have $$\mathbf{P}_x(\tau^+_y(X)< t)>0.$$
\edlemma
\proof
Since the case of $ \pi((0,\infty)) >0$ has been treated in Lemma~\ref{upup}, we may assume without loss of generality that $\pi((0,\infty)) =0$. Then the equation can be rewritten as
\beqnn
X_t= X_0+\int_0^t R_0(X_s)\,\d s+\int_0^t\sqrt{R_1(X_s)}\,\d B_s.
\eeqnn
For $0<\epsilon<x<\infty,$ denote $\hat\sigma:=\inf\{t\ge 0:X_t\notin[\epsilon/2,2y]\}.$ The scale function of $X_t$ is given by
$$\Phi(z):=\int_{\epsilon/2}^z\exp\left(-2\int_{\epsilon/2}^u\frac{ R_0(v)}{ R_1(v)}\d v\right)\d u,\quad z\in({\epsilon/2},\infty).$$
By It\^o's formula,
$$\Phi(X_{t\wedge \hat\sigma})=\Phi(x)+\int_0^{t\wedge \hat\sigma}\Phi'(X_s)\sqrt{ R_1(X_s)}\d B_s.$$
Let $M^\Phi_t=\int_0^{t\wedge \hat\sigma}\Phi'(X_s)\sqrt{ R_1(X_s)}\d B_s.$ Then $(M^\Phi_t)_{t\ge0}$ is a local martingale whose quadratic variation is
$$\langle M^\Phi\rangle_{t}=\int_0^{t\wedge \hat\sigma}\left[\Phi'(X_s)\sqrt{ R_1(X_s)}\right]^2\d s.$$
Define 
$$\theta(t)=\inf\{r\ge 0:\langle M^\Phi\rangle_{r}> t\},$$
and \beqlb\label{10010}\hat B_t=M^\Phi_{\theta(t)}=\Phi(X_{\theta(t)})-\Phi(x).\eeqlb
By  \cite[Theorem 4.6]{KS91}, $$\left\{\hat B_s,0<s<\langle M^\Phi\rangle_{\hat\sigma}\right\}$$ is a standard Brownian motion. Let
$$a=\Phi({\epsilon})-\Phi(x),\quad b=\Phi(y)-\Phi(x).$$
By \eqref{10010},
$$\tau_{b}^{a,+}(\hat B)=\int_0^{\tau_{b}^{a,+}(M^\Phi)}\left[\Phi'(X_s)\sqrt{ R_1(X_s)}\right]^2\d s\ge \hat C\tau_{b}^{a,+}(M^\Phi)$$
where $\hat C=\inf_{z\in(\epsilon,y)}\left[\Phi'(z)\sqrt{ R_1(z)}\right]^2.$
Since $\Phi$ is increasing on $(\epsilon/2,\infty),$
$$\mathbf P_x\left(\tau_{y}^{\epsilon,+}(X)<t\right)=\mathbf P_x\left(\tau_{b}^{a,+}(M^\Phi)<t\right)\ge\mathbf P_x\left(\tau_{b}^{a,+}(\hat B)<\hat Ct\right)>0.$$
\qed

\bgproposition[Upward reachability]\label{propupward}
Suppose that either \hyperlink{cond:C1}{\textup{(C1)}} holds or \(\pi((0,\infty))>0\). Then for any \(y>x>0\) and any \(t>0\),
\[
\mathbf P_x(\tau_y^+(X)<t)>0.
\]
\edproposition
\proof
If \(\pi((0,\infty))>0\), the assertion follows from Lemma \ref{upup}. If \(\pi((0,\infty))=0\), then the assertion follows from Lemma \ref{upupri} under \hyperlink{cond:C1}{\textup{(C1)}}.
\qed

\subsection{Downward reachability}

We now turn to downward reachability. This is the part where the time-change is needed: since the process has no negative jumps, a downward passage must be produced either by the continuous diffusion or by the compensator of truncated small jumps.

It remains to prove that for any $x>\epsilon>0$ and any $t_0>0,$  $\mathbf{P}_x(\tau^-_\epsilon(X)< t_0)>0$. The proof relies on the following time-change transformation.

Let $\bar R_i(x)=R_i(x)/R_2(x)$ for any $x\in(0,\infty)$ and $i=0,1$. 
For any $x,y\in(0,\infty),$
\beqnn
\left|\bar R_0(x)-\bar R_0(y)\right|\ar=\ar\left|\frac{1}{R_2(x)R_2(y)}\right|\left|R_0(x)R_2(y)-R_0(y)R_2(x)\right|\cr\cr
\ar\le\ar\left|\frac{1}{R_2(x)R_2(y)}\right|\left|R_0(x)R_2(x)-R_0(y)R_2(x)\right|\cr
\ar\ar+\left|\frac{1}{R_2(x)R_2(y)}\right|\left|R_0(x)R_2(y)-R_0(x)R_2(x)\right|\cr\cr
\ar=\ar\left|\frac{1}{R_2(y)}\right|\left|R_0(x)-R_0(y)\right|+\left|\frac{R_0(x)}{R_2(x)R_2(y)}\right|\left|R_2(y)-R_2(x)\right|.
\eeqnn
Thus for any closed interval $A\subset (0,\infty),$ there is a constant $C'(A)>0$ such that for any $x,y\in A,$
\beqlb\label{092801}|\bar R_0(x)-\bar R_0(y)|+|\bar R_1(x)-\bar R_1(y)|\le C'(A)|x-y|.\eeqlb
Let $(W_t)_{t\ge0}$ be a Brownian motion on $(\Omega,\mathcal F,\mathcal F_t,\mathbf P)$ and let $M(\d s,\d z,\d u)$ be a Poisson random measure on the same probability space with intensity $\d s\pi(\d z)\d u.$ Then by \cite[Theorem 3.1]{LYZ}, there is a unique strong solution $(Z'_t)_{t\ge0}$ of the following SDE, up to the first time of hitting $0$ or explosion. 
\beqlb\label{092001}
Z'_{t}\ar=\ar x+\int_0^{t}\bar R_0(Z'_s)\d s+\int_{0}^{t} \sqrt{\bar R_{1}\left(Z'_{s}\right)} \d W_{s} +\int_0^{t}\int_0^\infty\int_0^1 z \tilde{M}(\d s,\d z,\d u)
\eeqlb

The purpose of this time change is to remove the state dependence from the jump intensity. After this time change, the small-jump part and the large-jump part of the Poisson noise are independent, which allows us below to construct truncated auxiliary paths that coincide with the original path with positive probability.

\bglemma\label{translem0923}
Let $Z_t=Z'_{t\wedge\tau_0^-(Z')\wedge\tau_\infty^+(Z')},$ and
\[
\kappa(t)=\inf\left\{u\ge 0:\int_0^u \frac{1}{R_2(Z_s)}\d s>t\right\}.
\]
Then \(X_t:=Z_{\kappa(t)}\) solves the equation \eqref{eq:sde_intro} with
\(X_0=x\) on an extension of the original probability space, up to its
lifetime. By the standing pathwise uniqueness assumption for
\eqref{eq:sde_intro}, this time-changed solution has the same law as the
original solution started from \(x\). In particular, the identity holds before every hitting time
\(\tau^-_\epsilon(X)\), \(\epsilon>0\).
\edlemma
\begin{proof}
Let \(X_t=Z_{\kappa(t)}\). For \(n\ge1\), set
\[
\zeta_n=\inf\left\{t\ge 0: X_t\ge n \text{ or } X_t\le \frac{1}{n}\right\},
\qquad
\zeta=\lim_{n\to\infty}\zeta_n.
\]
We work on \([0,\zeta_n]\), where \(R_2(X)\) is bounded away from zero.
Since the clock \(A_u:=\int_0^u
R_2(Z_s)^{-1}\d s\) has an inverse \(\kappa\), we have
\[
\kappa(t\wedge\zeta_n)=\int_0^{t\wedge\zeta_n}R_2(X_{s-})\d s.
\]
Therefore,
{\small\beqlb
X_{t\wedge\zeta_n} \ar=\ar x
+ \int_0^{\kappa(t\wedge\zeta_n)} \bar R_0(Z_s)\d s
+ \int_0^{\kappa(t\wedge\zeta_n)} \sqrt{\bar R_1(Z_s)}\d W_s \ccr
\ar\ar\quad
+ \int_0^{\kappa(t\wedge\zeta_n)}\int_0^\infty\int_0^1 z\,\tilde M(\d s,\d z,\d u).
\label{trans1}
\eeqlb}
Moreover, by the clock identity and the definition \(\bar R_0=R_0/R_2\),
\beqlb\label{drift0923}
\int_0^{\kappa(t\wedge\zeta_n)} \bar R_0(Z_s)\d s
=\int_0^{t\wedge\zeta_n}R_0(X_s)\d s.
\eeqlb

We next make the Brownian motion on the new time scale explicit. For
\(t\ge0\), define the stopped local martingale
\[
B_{t\wedge\zeta_n}
:=\int_0^{\kappa(t\wedge\zeta_n)}
\frac{1}{\sqrt{R_2(Z_s)}}\,\d W_s .
\]
Its quadratic variation is
\[
\left\langle B_{\cdot\wedge\zeta_n}\right\rangle_t
=\int_0^{\kappa(t\wedge\zeta_n)}\frac{1}{R_2(Z_s)}\,\d s
=t\wedge\zeta_n.
\]
Hence, by L\'evy's martingale characterization (see
\cite[Chapter II, Theorem 6.1]{ikeda}), there exists a Brownian motion \(B\)
on an extension of the probability space whose stopped versions are the stopped
local martingales defined above. The stochastic
change-of-variables formula then gives
\beqlb\label{W}
\int_0^{\kappa(t\wedge\zeta_n)} \sqrt{\bar R_1(Z_s)}\,\d W_s
=\int_0^{t\wedge\zeta_n}\sqrt{R_1(X_s)}\,\d B_s .
\eeqlb

We now transform the Poisson random measure. Take an additional Poisson
random measure \(M_1(\d s,\d z,\d u)\), independent of the original noises,
with intensity \(\d s\,\pi(\d z)\,\d u\), and define \(N\) by
\beqlb\label{N1}
N(\d s,\d z,\d u):\ar=\ar
\I_{\{s<\zeta,\ u\le R_2(X_{s-})\}}
M\bigl(\d\kappa(s),\d z,[R_2(X_{s-})]^{-1}\d u\bigr)\cr
\ar\ar\qquad\qquad
+\I_{\{s\ge\zeta\}\cup\{u>R_2(X_{s-})\}}M_1(\d s,\d z,\d u).
\eeqlb
The martingale characterization of Poisson point
processes \cite[Chapter II, Theorem 6.2]{ikeda} shows that \(N\) is a
Poisson random measure with intensity \(\d s\,\pi(\d z)\,\d u\). Moreover,
applying the joint characterization in
\cite[Chapter II, Theorem 6.3]{ikeda}, the Brownian motion \(B\) and the Poisson random
measure \(N\) are independent. Finally,
\[
\int_0^{\kappa(t\wedge\zeta_n)}\int_0^\infty\int_0^1
z\,\tilde M(\d s,\d z,\d u)
=
\int_0^{t\wedge\zeta_n}\int_0^\infty\int_0^{R_2(X_{s-})}
z\,\tilde N(\d s,\d z,\d u).
\]
Combining this identity with \eqref{trans1}, \eqref{drift0923}, and \eqref{W},
we obtain
\[
X_{t\wedge\zeta_n}
=x+\int_0^{t\wedge\zeta_n}R_0(X_s)\d s
+\int_0^{t\wedge\zeta_n}\sqrt{R_1(X_s)}\,\d B_s
+\int_0^{t\wedge\zeta_n}\int_0^\infty\int_0^{R_2(X_{s-})}
z\,\tilde N(\d s,\d z,\d u).
\]
Letting \(n\to\infty\) proves the desired time-change representation up to
the lifetime, and in particular on every interval before \(\tau^-_\epsilon(X)\).
\qed
\end{proof}
\bglemma\label{lem0923}
Suppose that $\int_0^1 z\pi(\d z)=\infty$. Then for any $x>\epsilon>0$ and any $T>0,$ we have $\mathbf{P}_x(\tau^-_\epsilon(Z)< T)>0$.
\edlemma

\proof 
Fix \(T>0\). Choose \(t_0\in(0,T)\), to be taken sufficiently small below. It is enough to prove \(\mathbf{P}_x(\tau^-_\epsilon(Z)< t_0)>0\).
For $\delta\in (0, 1)$, we introduce an auxiliary process defined by the following equation: 
\beqlb\label{auxil1} 
Y^\delta_{t}= x+\int_0^{t}\bar{R}_0(Y^\delta_s)\d s+M^\delta_t-\int^t_0 \left( \int_\delta^\infty z \pi(\d z) \right)\d s,
\eeqlb
where $(M^\delta_t)_{t\ge0}$ is a local martingale given by
\beqnn
M^\delta_t:=\int_{0}^{t} \sqrt{\bar{R}_{1}\left(Y^\delta_{s}\right)} \d W_{s} +\int_0^{t}\int_0^\delta\int_0^1 z \tilde{M}(\d s,\d z,\d u).
\eeqnn
We shall show that for any fixed short time interval $(0,t_0)$, one can choose $\delta\in(0,1)$ such that the auxiliary process $(Y^\delta_t)_{t\ge0}$ reaches the level \(\epsilon\) and pathwise coincides with the original process on \((0,t_0)\) with positive probability.

We first prove that, for this sufficiently small \(t_0\),
\beqlb\label{092301}
\mathbf P_x\left(\tau^-_\epsilon(Y^\delta)< t_0\right)>0,
\eeqlb
where $\delta$ will be specified later. To establish this, we consider the stopped process $(Y^\delta_{t\wedge\sigma^\delta})_{t\ge0}$ where the stopping time is defined as
\beqnn
\sigma^\delta:=\tau^-_{{\epsilon}/{2}}(Y^\delta) \wedge \tau_{x+2}^+(Y^\delta).
\eeqnn
By the definition in \eqref{auxil1}, it suffices to prove that the sum of the probabilities of the following three events is strictly less than $1$:
\beqlb\label{9271}
1\ar>\ar
\mathbf P_x\bigg(\int_0^{t_0\wedge\sigma^\delta} \left|\bar{R}_0(Y^\delta_s)\right|\d s\ge 1\bigg)+\mathbf P_x\left( \left|M_{t_0\wedge\sigma^\delta}^\delta\right|\ge 1\right)\cr
\ar\ar+\mathbf{P}_x\bigg(\int_0^{t_0\wedge\sigma^\delta}\left(\int^\infty_\delta z \pi(\d z)\right)\d s<x-\epsilon+2, \;\; \tau^-_{{\epsilon}/{2}}(Y^\delta)>{t_0\wedge\sigma^\delta}\bigg)\cr
\ar=:\ar p_1+p_2+p_3.
\eeqlb
We now proceed to estimate $p_1, p_2$ and $p_3$ separately. 

We begin with $p_1$. Since $Y^\delta_s \in [\epsilon/2, x+2]$ up to $\sigma^\delta$, $\bar{R}_0$ is uniformly bounded on this compact interval. Thus, we can choose a sufficiently small $t_0$ such that the integral is deterministically less than $1$, yielding
\beqlb\label{092207}
p_1\le\mathbf{P}_x\left(\int_0^{t_0\wedge\sigma^\delta}\left|\bar{R}_0(Y^\delta_s)\right|\d s\ge1\right)=0.
\eeqlb

Next, we turn to estimate $p_2.$
Note that $(M_{t\wedge\sigma^\delta}^\delta)_{t\ge0}$ is a martingale whose predictable quadratic variation satisfies
\beqlb\label{092091}
\langle M^\delta\rangle_{t\wedge\sigma^\delta} \le \int_0^{t\wedge\sigma^\delta}\bar{R}_1(Y^\delta_s)\d s+(t\wedge\sigma^\delta)\int_0^\delta z^2\pi(\d z).
\eeqlb
By applying Chebyshev's inequality and It\^o's isometry, we have
\beqnn
p_2\le\mathbf{P}_x\left(\left|M_{t_0\wedge\sigma^\delta}^\delta\right|\ge 1\right) \ar\le\ar \mathbf E_x\left[\left|M_{t_0\wedge\sigma^\delta}^\delta\right|^2\right]\le\left(\sup_{y \in [\epsilon/2, x+2]}\bar{R}_1(y)+\int_0^{1}z^2\pi(\d z)\right) t_0.
\eeqnn
Thus, by further reducing $t_0$ if necessary, we obtain
\beqlb\label{092205}
p_2 \le \frac{1}{3}.
\eeqlb

Finally, we estimate $p_3.$
The infinite variation condition $\int^1_0 z \pi(\d z)=\infty$ implies that we can find a $\delta\in (0, 1)$ satisfying
\beqlb\label{091005}
\int_\delta^\infty z\pi(\d z) \ge \frac{x-\epsilon+2}{t_0}.
\eeqlb
Then, on the event $\{\sigma^\delta\ge t_0\}$, we deterministically have 
\beqnn
\int_0^{t_0\wedge\sigma^\delta}\left(\int^\infty_\delta z \pi(\d z)\right)\d s \ge x-\epsilon+2.
\eeqnn
Therefore, the event in $p_3$ implies that $\sigma^\delta < t_0$. Furthermore, since $\tau^-_{{\epsilon}/{2}}(Y^\delta)>\sigma^\delta$, it must be that $\sigma^\delta = \tau_{x+2}^+(Y^\delta)$, which means $Y^\delta_{\sigma^\delta} - x \ge 2$. Thus,
\beqlb\label{092302}
p_3 \ar\le\ar \mathbf P_x( \sigma^\delta < t_0, \;\; Y^\delta_{\sigma^\delta}-x\ge 2)\cr
\ar\le\ar\mathbf P_x\left(\int_0^{t_0\wedge\sigma^\delta}\left|\bar{R}_0(Y^\delta_s)\right|\d s + M^\delta_{t_0\wedge\sigma^\delta} \ge 2\right)\cr
\ar\le\ar\mathbf P_x\left(\int_0^{t_0\wedge\sigma^\delta}\left|\bar{R}_0(Y^\delta_s)\right|\d s\ge1\right)+\mathbf P_x\left(\left|M^\delta_{t_0\wedge\sigma^\delta}\right|\ge1\right) \cr
\ar\le\ar p_1 + p_2.
\eeqlb
Substituting \eqref{092207} and \eqref{092205} into \eqref{092302}, we have
\beqlb\label{0922v1}
p_3\le 0 + \frac{1}{3} = \frac{1}{3}.
\eeqlb
Combining \eqref{092207}, \eqref{092205} and \eqref{0922v1} gives $p_1 + p_2 + p_3 \le 2/3 < 1$, which implies 
$$
\mathbf P_x(\tau^-_\epsilon(Y^\delta)< t_0)>0.
$$

Having established that the auxiliary process can reach level $\epsilon$ within the short time interval $(0, t_0]$, we now prove its pathwise coincidence with the original process. Rewrite equation \eqref{092001} as 
\beqnn
Z_{t}\ar=\ar x+\int_0^{t}\bar R_0(Z_s)\d s+\int_{0}^{t} \sqrt{\bar R_{1}\left(Z_{s}\right)} \d W_{s} +\int_0^{t}\int_0^\delta\int_0^1 z \tilde{M}(\d s,\d z,\d u)\\
\ar\ar-\int^t_0\int_\delta^\infty z \pi(\d z)\d s+\int_0^{t}\int_\delta^\infty\int_0^1 z M(\d s,\d z,\d u).
\eeqnn
Let $\sigma_1=\inf\{t\ge 0: \int_0^{t}\int_{\delta}^\infty\int_0^1 M(\d s,\d z,\d u)>0\}$. Then for any $t_0\in(0,\infty),$
\beqnn
\mathbf{P}_x(\sigma_1>t_0)=\mrm{e}^{- \pi((\delta, \infty))t_0}>0.
\eeqnn
By the spatial independence of Poisson random measures (e.g., \cite[Chapter II, Theorem 6.3]{ikeda}), the jump mechanisms inside and outside the region $z \in (0, \delta]$ are independent. Thus, $(Y^\delta_t)_{t\ge0}$ and $\sigma_1$ are independent. On the event \(\{\sigma_1>t_0\}\), the processes \(Z\) and \(Y^\delta\) solve the same truncated equation on \([0,t_0]\); by pathwise uniqueness, \(Z_t=Y^\delta_t\) for \(t\le t_0\). As a consequence,
$$
\mathbf{P}_x\left(\tau^-_\epsilon(Z)<t_0\right)\ge \mathbf{P}_x\left(\tau^-_\epsilon(Y^\delta)<t_0, \;\; \sigma_1>t_0\right) = \mathbf{P}_x\left(\tau^-_\epsilon(Y^\delta)<t_0\right) \mathbf{P}_x\left(\sigma_1>t_0\right)>0.
$$ 
This completes the proof.
\qed

\bglemma\label{lemdown}
Suppose that $R_1$ is strictly positive on $(0,\infty)$ and \(\int_0^1 z\pi(\d z)<\infty\). Then for any $T>0$ and any $0<\epsilon<x,$ $\mathbf P_x(\tau_\epsilon^-(Z)<T)>0.$
\edlemma
\proof
We first introduce the continuous diffusion \(Y\) by
\beqlb\label{auxil2loc}
Y_t=x+\int_0^t\bigl[\bar R_0(Y_s)-m_\pi\bigr]\d s+\int_0^t\sqrt{\bar R_1(Y_s)}\,\d W_s,
\eeqlb
where \[
m_\pi:=\int_0^\infty z\pi(\d z)<\infty.
\]
For \(\delta\in(0,1)\), define the truncated process \(Z^\delta\) by
\beqlb\label{auxil3loc}
Z^\delta_t\ar=\ar x+\int_0^t\bigl[\bar R_0(Z^\delta_s)-m_\pi\bigr]\d s+\int_0^t\sqrt{\bar R_1(Z^\delta_s)}\,\d W_s
+\int_0^t\int_0^\delta\int_0^1 z M(\d s,\d z,\d u).
\eeqlb
Since \(\bar R_1\) is strictly positive on compact subintervals of \((0,\infty)\), the local irreducibility of one-dimensional non-degenerate diffusions gives, for some \(M>x\vee1\),
\beqlb\label{eq:loc-p0}
p_0:=\mathbf P_x(E)>0,\qquad
E:=\left\{\tau^-_{\epsilon/2}(Y)<T\wedge\tau^+_M(Y)\right\}.
\eeqlb
Set
\[
\rho^\delta:=\tau^-_{\epsilon/2}(Y)\wedge\tau^+_M(Y)\wedge\tau^+_{M+\epsilon/2}(Z^\delta)\wedge T.
\]
The comparison theorem \cite[Theorem 5.5]{FL10} gives \(Y_t\le Z^\delta_t\) for \(0\le t\le\rho^\delta\).
On \([0,\rho^\delta)\), the two processes are localized below \(M+\epsilon/2\); moreover, the jumps of \(Z^\delta\) are bounded by \(1\). Hence, by \eqref{092801}, there is a constant \(C_M>0\) such that for all \(u,v\in[\epsilon/2,M+\epsilon/2+1]\),
\[
|\bar R_0(u)-\bar R_0(v)|+\left|\sqrt{\bar R_1(u)}-\sqrt{\bar R_1(v)}\right|^2\le C_M|u-v|.
\]
Let \(\eta^\delta_t=Z^\delta_t-Y_t\ge0\). Subtracting \eqref{auxil2loc} from \eqref{auxil3loc}, stopping at \(t\wedge\rho^\delta\), and taking expectations yield
\beqlb\label{eq:loc-eta}
\mathbf E_x\left[\eta^\delta_{t\wedge\rho^\delta}\right]
\le C_M\int_0^t\mathbf E_x\left[\eta^\delta_{s\wedge\rho^\delta}\right]\d s
+T\int_0^\delta z\pi(\d z),\qquad 0\le t\le T.
\eeqlb
By Gronwall's inequality,
\beqlb\label{eq:loc-gronwall}
\mathbf E_x\left[\eta^\delta_{\rho^\delta}\right]
\le T \mrm e^{C_MT}\int_0^\delta z\pi(\d z).
\eeqlb
Choosing \(\delta>0\) sufficiently small, we may assume
\[
\mathbf E_x\left[\eta^\delta_{\rho^\delta}\right]
<\frac{\epsilon p_0}{4}.
\]
Let \(G_\delta:=\left\{\eta^\delta_{\rho^\delta}<{\epsilon}/{2}\right\}.\) By Markov's inequality and the choice of \(\delta\),
\[
\mathbf P_x(G_\delta^c)\le
\frac{2}{\epsilon}
\mathbf E_x\left[\eta^\delta_{\rho^\delta}\right]
<\frac{p_0}{2}.
\]
By \eqref{eq:loc-p0},
\[
\mathbf P_x(E\cap G_\delta)
\ge \mathbf P_x(E)-\mathbf P_x(G_\delta^c)
>\frac{p_0}{2}>0.
\]
On \(E\cap G_\delta\), the process \(Z^\delta\) cannot reach \(M+\epsilon/2\) before \(\tau^-_{\epsilon/2}(Y)\). Indeed, if this happened, then at \(\rho^\delta\) we would have \(Z^\delta_{\rho^\delta}\ge M+\epsilon/2\) while \(Y_{\rho^\delta}\le M\), and hence \(\eta^\delta_{\rho^\delta}\ge\epsilon/2\), contradicting \(G_\delta\). Thus, on \(E\cap G_\delta\),
\[
\rho^\delta=\tau^-_{\epsilon/2}(Y)<T.
\]
Since \(Y\) is continuous, \(Y_{\rho^\delta}=\epsilon/2\). Together with \(G_\delta\), this gives
\[
Z^\delta_{\rho^\delta}
=Y_{\rho^\delta}+\eta^\delta_{\rho^\delta}
<\epsilon.
\]
Since \(Z^\delta_0=x>\epsilon\) and \(Z^\delta\) has no negative jumps, this implies \(E\cap G_\delta\subset\{\tau^-_\epsilon(Z^\delta)<T\}\). Consequently,
\[
\mathbf P_x\left(\tau^-_\epsilon(Z^\delta)<T\right)
\ge \mathbf P_x(E\cap G_\delta)>0.
\]

We now compare \(Z^\delta\) with the original time-changed process \(Z\). In the present finite-variation case, \eqref{092001} can be written as
\beqlb\label{eq:Z-fv}
Z_t\ar=\ar x+\int_0^t\bigl[\bar R_0(Z_s)-m_\pi\bigr]\d s+\int_0^t\sqrt{\bar R_1(Z_s)}\,\d W_s\cr
\ar\ar+\int_0^t\int_0^\delta\int_0^1 z M(\d s,\d z,\d u)
+\int_0^t\int_\delta^\infty\int_0^1 z M(\d s,\d z,\d u).
\eeqlb
Thus the equation for \(Z\) differs from \eqref{auxil3loc} only by the last large-jump term. Let
\[
\sigma^\delta=\inf\left\{t\ge0:\int_0^t\int_\delta^\infty\int_0^1 M(\d s,\d z,\d u)>0\right\}.
\]
Since \(\pi((\delta,\infty))<\infty\), \(\mathbf P_x(\sigma^\delta>T)=\exp\{-\pi((\delta,\infty))T\}>0\). The event \(\{\sigma^\delta>T\}\) is independent of the noises appearing in \eqref{auxil3loc}. On this event, the last term in \eqref{eq:Z-fv} is zero on \([0,T]\), and hence \(Z\) and \(Z^\delta\) coincide up to \(T\) by pathwise uniqueness. Thus
\[
\mathbf P_x(\tau^-_\epsilon(Z)<T)
\ge \mathbf P_x(\tau^-_\epsilon(Z^\delta)<T,\sigma^\delta>T)
= \mathbf P_x(\tau^-_\epsilon(Z^\delta)<T)\mathbf P_x(\sigma^\delta>T)>0.
\]
\qed

\bgproposition[Downward reachability]\label{propddd}
Suppose that either \hyperlink{cond:C1}{\textup{(C1)}} or \hyperlink{cond:C3}{\textup{(C3)}} holds. Then for any $x>\epsilon>0$ and any $t_0>0,$ we have $\mathbf{P}_x(\tau^-_\epsilon(X)< t_0)>0$.
\edproposition
\proof
Recall the time-change relationship between $X$ and the auxiliary process $Z$. By the monotonicity of $R_2$, we have $R_2(Z_s) \ge R_2(\epsilon)$ for all $0 \le s \le \tau_\epsilon^-(Z)$. Therefore,
$$\tau_\epsilon^{-}(X)=\int_0^{\tau_\epsilon^-(Z)}\frac{1}{R_2(Z_s)}\d s\le \frac{\tau_\epsilon^-(Z)}{R_2(\epsilon)}.$$
Thus, for any $t_0 > 0,$
$$\mathbf{P}_x(\tau^-_\epsilon(X)< t_0)\ge \mathbf{P}_x(\tau^-_\epsilon(Z)< t_0 R_2(\epsilon)).$$
If \hyperlink{cond:C3}{\textup{(C3)}} holds, the strict positivity of the right-hand side follows from Lemma \ref{lem0923}. If \hyperlink{cond:C1}{\textup{(C1)}} holds and \hyperlink{cond:C3}{\textup{(C3)}} fails, then \(\int_0^1 z\pi(\d z)<\infty\), and Lemma \ref{lemdown} applies because \(R_1\) is strictly positive on \((0,\infty)\).
\qed

{\it Proof of Theorem \ref{irred}}
Under the present assumptions, Propositions \ref{propupward} and \ref{propddd} verify the two one-sided hitting hypotheses of Proposition \ref{prop_reach_irred}. Hence \(P_t\) is irreducible on \((0,\infty)\).
\qed

\bgproposition[Interior regularity package]\label{prop_regular_package}
Suppose that either \hyperlink{cond:C1}{\textup{(C1)}} holds, or both
\hyperlink{cond:C2}{\textup{(C2)}} and \hyperlink{cond:C3}{\textup{(C3)}} hold. Then the original process \(X\) has the following properties.
\begin{itemize}
    \item For any \(t>0\) and any \(x_0\in(0,\infty)\),
    \[
    \lim_{x\to x_0}\|P_t(x,\cdot)-P_t(x_0,\cdot)\|_{\mathrm{TV}}=0.
    \]
    \item For any \(T>0\), the map
    \[
    x\mapsto \mathbf P_x\bigl((X_t)_{t\in[0,T]}\in\cdot\bigr)
    \]
    is continuous on \((0,\infty)\) in the sense of weak convergence on \(D([0,T],\mathbb R_+)\).
    \item For any compact set \(K\subset(0,\infty)\) and any \(T>0\), with
    \(\tau_R^+=\inf\{t\ge0:X_t\ge R\}\),
    \[
    \lim_{R\to\infty}\sup_{x\in K}\mathbf P_x(\tau_R^+\le T)=0.
    \]
    \item For any \(x>0\), \(t>0\), and nonempty open set \(O\subset(0,\infty)\),
    \[
    P_t(x,O)>0.
    \]
\end{itemize}
\edproposition
\proof
Under the stated alternatives, the total-variation Feller property follows from Theorem \ref{SF_a}, the trajectory Feller property from Theorem \ref{trf}, the compact-starting finite-time boundedness from Proposition \ref{prop_compact_bounded}, and the open-set support positivity from Theorem \ref{irred}.
\qed

\section{Applications}

\subsection{Quasi-stationary distribution}
We first apply the strong Feller property and irreducibility to the quasi-stationary distribution (QSD) problem. Together with a Lyapunov condition, these topological properties yield existence, uniqueness, and exponential convergence of the QSD; see \cite[Theorem 2.2]{GNW20}. We begin by verifying the required Lyapunov condition.

\bgproposition\label{lyp}
Suppose that Condition \hyperlink{cond:C4}{\textup{(C4)}} holds. Then there exists a function $\psi\in C_b^2([0,\infty))$ with $\psi(x)\ge 1$ for all $x\in[0,\infty),$ sequences of positive constants $(r_n)_{n\ge 1}$ and $(b_n)_{n\ge 1},$ with $r_n\rightarrow\infty$ as $n\rightarrow\infty,$ and an increasing sequence $(K_n)_{n\ge 1}$ of compact subsets of $[0,\infty),$ such that for all $n\ge 1,$
$$L\psi(x)\le -r_n\psi(x)+b_n\I_{K_n}(x),\quad x\in[0,\infty).$$
\edproposition
\begin{proof}
Define the function $\psi\in C_b^2([0, \infty))$ by
\beqnn
\psi(x) = 2 - (x+1)^{-r/2}, \quad x \ge 0,
\eeqnn
where $r > 0$ is the constant given in Condition \hyperlink{cond:C4}{\textup{(C4)}}. Notice that for all $x \ge 0,$
\beqnn
1\le \psi(x)<2,\qquad
\psi'(x)=\frac r2(x+1)^{-\frac r2-1}>0,\qquad
\psi''(x)=-\frac r2\left(\frac r2+1\right)(x+1)^{-\frac r2-2}<0.
\eeqnn
Applying the generator $L$ defined by \eqref{gen} to $\psi$, we have
\beqlb\label{eq:psinbd}
L\psi(x) \ar=\ar R_0(x)\psi'(x) + \frac{1}{2}R_1(x)\psi''(x) \cr\ar\ar+ R_2(x) \int_0^\infty \left[ \psi(x+z) - \psi(x) - z\psi'(x) \right] \pi(\d z)\cr
\ar\le\ar
R_0(x)\psi'(x).
\eeqlb
By Condition \hyperlink{cond:C4}{\textup{(C4)}}, we have
\beqlb\label{eqpsir0}
R_0(x)\psi'(x)
\ar\le\ar \frac r2\left(-\lambda_1 x^{r+1} + \lambda_2\right)(x+1)^{-\frac r2-1}\cr
\ar\le\ar -\frac r2\lambda_1 x^{r+1}(2x)^{-\frac r2-1}
+\frac r2\lambda_2\cr
\ar=\ar -2^{-\frac r2-2}r\lambda_1 x^{r/2}
+\frac r2\lambda_2,
\qquad x\ge1.
\eeqlb
Thus, taking $c=2^{-\frac r2-3}r\lambda_1,$ there exists a sufficiently large $l>0$  such that 
\[L\psi(x)\le-c x^{r/2},\qquad x\ge l.\]
Since $1\le \psi(x)<2,$ for $x\ge l,$
\[L\psi(x)\le -cx^{r/2}\le-\frac c2x^{r/2}\psi(x).\]

For $n=1,2,\cdots,$ define
\[K_n=[0,l+n-1],\quad r_n=\frac c2(l+n-1)^{r/2},\quad b_n=2r_n+\sup_{x\in K_n}|L\psi(x)|.\]
For $x\notin K_n,$
\[L\psi(x)\le-\frac c2x^{r/2}\psi(x)\le- r_n\psi(x).\]
For $x\in K_n,$ since $\psi(x)<2,$
\[-r_n\psi(x)+b_n\ge -2r_n+(2r_n+\sup_{x\in K_n}|L\psi(x)|)=\sup_{x\in K_n}|L\psi(x)|\ge L\psi(x).\]
This completes the proof.

\qed
 \end{proof}
{\it Proof of Theorem \ref{qsd}.} Under \hyperlink{cond:C5}{\textup{(C5)}}, by \eqref{es04292}, for any fixed $t>0,$ we have
$$
\mathbf{P}_x \left(\tau^-_0>t\right)\to0\quad\mbox{as}\quad x\to 0.
$$
In particular, there exists \(x_0>0\) such that
$$\mathbf{P}_{x_0}(\tau^-_0<\infty)>0,$$
which verifies the non-trivial killing condition in \cite[Theorem 2.2]{GNW20}.
Under the stated alternatives, Theorems \ref{SF_a} and \ref{irred} give the strong Feller property and irreducibility in \((0,\infty)\).
Let \(D=(0,\infty)\), and denote the killed semigroup by
\[
P_t^D f(x):=\mathbf E_x\left[f(X_t)\I_{\{t<\tau_0^-\}}\right],
\qquad x\in D.
\]
For \(f\in b\mathcal B(D)\), extend \(f\) to \([0,\infty)\) by setting \(\tilde f(0)=0\) and \(\tilde f=f\) on \(D\). Since \(0\) is absorbing, \(P_t^D f(x)=P_t\tilde f(x)\) for \(x\in D\). Therefore the strong Feller and boundary-continuity inputs for the killed semigroup follow from Theorems \ref{SF_a} and \ref{sf0}. Moreover, for every nonempty open set \(O\subset D\),
\[
P_t^D(x,O)=P_t(x,O),\qquad x\in D,\ t>0,
\]
and hence the killed irreducibility required in \cite[Theorem 2.2]{GNW20} follows from Theorem \ref{irred}. The remaining conditions required by \cite[Theorem 2.2]{GNW20} are checked as follows: Theorems \ref{trf} and \ref{trf0} give the trajectory Feller condition, and Proposition \ref{lyp} verifies the Lyapunov condition. Hence, all hypotheses of \cite[Theorem 2.2]{GNW20} are satisfied, and Theorem \ref{qsd} follows.

\qed

\subsection{Uniform exponential ergodicity}

The following boundary lemma will be used to handle compact sets touching the boundary. We include the proof for completeness.
\bglemma\label{lemzero}
Suppose that \(R_0(0)>0\) and that \((P_t)_{t\ge0}\) is irreducible on \((0,\infty)\). Then the following assertions hold.
\begin{enumerate}[label=\rm(\roman*)]
\item For any \(t>0\) and any nonempty open set \(O\subset(0,\infty)\), we have \(P_t(0,O)>0\).
\item If, in addition, \((P_t)_{t\ge0}\) is strong Feller on \((0,\infty)\), then there exist \(t_0>0\), \(\delta>0\), \(c_0>0\), and a compact interval \(K_0\subset(0,\infty)\) with nonempty interior such that
\[
\inf_{y\in[0,\delta]}P_{t_0}(y,K_0)\ge c_0.
\]
\end{enumerate}
\edlemma
\proof
{\rm(i).}
Fix \(t>0\) and a nonempty open set \(O\subset(0,\infty)\), and set \(s=t/2\). We first show that the process can leave the boundary with positive probability before time \(s\). For \(\lambda>0\), set \(g_\lambda(x)=\exp(-\lambda x)\). Then for \(x>0\),
\[
Lg_\lambda(x)=g_\lambda(x)\left[-\lambda R_0(x)+\frac{\lambda^2}{2}R_1(x)+R_2(x)J_\lambda\right],
\]
where\[
J_\lambda:=\int_0^\infty \left(\mrm e^{-\lambda z}-1+\lambda z\right)\pi(\d z)<\infty.
\]
Since \(R_1(0)=R_2(0)=0\), the right-hand side extends continuously to \(x=0\), with value \(-\lambda R_0(0)\).
Choose \(c>1/s\). Since \(R_0(0)>0\), we may first choose \(\lambda\) sufficiently large and then \(\varepsilon>0\) sufficiently small so that, for all \(x\in[0,\varepsilon]\),
\[
Lg_\lambda(x)\le -c.
\]
Let \(\tau_\varepsilon^+=\inf\{r\ge0:X_r\ge\varepsilon\}\). For any \(y\in[0,\varepsilon]\), Dynkin's formula applied to \(g_\lambda\) up to \(s\wedge\tau_\varepsilon^+\) gives
\[
\mathbf E_yg_\lambda(X_{s\wedge\tau_\varepsilon^+})
=g_\lambda(y)+\mathbf E_y\int_0^{s\wedge\tau_\varepsilon^+}Lg_\lambda(X_r)\d r
\le g_\lambda(y)-c\mathbf E_y[s\wedge\tau_\varepsilon^+].
\]
Since \(g_\lambda\ge0\), it follows that
\beqlb\label{eq:lemzero-boundary-dynkin}
\mathbf E_y[s\wedge\tau_\varepsilon^+]\le \frac{g_\lambda(y)}{c}\le c^{-1},
\qquad y\in[0,\varepsilon].
\eeqlb
Consequently,
\[
\mathbf P_0(\tau_\varepsilon^+<s)\ge 1-\frac{1}{cs}>0.
\]
Applying the strong Markov property at \(\tau_\varepsilon^+\), we obtain
\[
P_t(0,O)\ge
\mathbf E_0\left[\I_{\{\tau_\varepsilon^+<t/2\}}
P_{t-\tau_\varepsilon^+}(X_{\tau_\varepsilon^+},O)\right]>0,
\]
because \(X_{\tau_\varepsilon^+}>0\) on \(\{\tau_\varepsilon^+<t/2\}\) and the assumed irreducibility gives
\(P_{t-\tau_\varepsilon^+}(X_{\tau_\varepsilon^+},O)>0\).

{\rm(ii).}
Fix \(s>0\), take \(c>4/s\), and choose \(\lambda>0\) and \(\varepsilon>0\) as
in {\rm(i)} so that \(Lg_\lambda(x)\le -c\) for \(x\in[0,\varepsilon]\). Then
\eqref{eq:lemzero-boundary-dynkin} gives
\[
\inf_{y\in[0,\varepsilon/2]}\mathbf P_y(\tau_\varepsilon^+<s)
\ge q:=1-\frac{1}{cs}>0.
\]
Moreover, applying the stopped equation up to \(s\wedge\tau_\varepsilon^+\) and
using that the compensated integral has mean zero, we get
\[
\mathbf E_y\left[X_{s\wedge\tau_\varepsilon^+}\right]
=y+\mathbf E_y\int_0^{s\wedge\tau_\varepsilon^+}R_0(X_r)\d r
\le \frac{\varepsilon}{2}+s\sup_{0\le x\le\varepsilon}|R_0(x)|=:B,
\qquad y\in[0,\varepsilon/2].
\]
Choose \(M>2\varepsilon\) so large that \(2B/M\le q/2\). By Markov's inequality,
\[
\inf_{y\in[0,\varepsilon/2]}
\mathbf P_y\left(\tau_\varepsilon^+<s,\ X_{\tau_\varepsilon^+}\le M/2\right)
\ge q/2>0.
\]

We next prove a uniform interior lower bound over the compact time interval
\([s,2s]\). We claim that there exists \(\beta>0\) such that
\beqlb\label{eq:lemzero-uniform-interior}
\inf_{x\in[\varepsilon,M/2]}\inf_{r\in[s,2s]}P_r(x,[\varepsilon/4,M])\ge \beta.
\eeqlb
Indeed, fix \(r_*\in[s,2s]\). By Lemma~\ref{lemineq2}, with
\[
(a,b)=(\varepsilon/4,M),\qquad (a',b')=(\varepsilon/2,3M/4),
\]
we may choose \(h\in(0,r_*/2)\) small enough so that
\[
\inf_{z\in[\varepsilon,M/2]}\inf_{0\le u\le 2h}P_u(z,[\varepsilon/4,M])\ge \frac12.
\]
Since \(P_{r_*-h}(x,(\varepsilon,M/2))\) is strictly positive by irreducibility and is
continuous in \(x\in(0,\infty)\) by the interior strong Feller property,
\[
\alpha_{r_*}:=\inf_{x\in[\varepsilon,M/2]}P_{r_*-h}(x,(\varepsilon,M/2))>0.
\]
Thus, for all \(r\in[r_*-h,r_*+h]\) and \(x\in[\varepsilon,M/2]\),
\[
\begin{split}
P_r(x,[\varepsilon/4,M])
&\ge \int_{(\varepsilon,M/2)}P_{r_*-h}(x,\d z)P_{r-(r_*-h)}(z,[\varepsilon/4,M])\\
&\ge \frac12\,\alpha_{r_*}.
\end{split}
\]
For each \(r_*\in[s,2s]\), let \(h(r_*)\in(0,r_*/2)\) be chosen as above.
Then the open intervals
\[
I_{r_*}:=(r_*-h(r_*),r_*+h(r_*)),\qquad r_*\in[s,2s],
\]
form an open cover of the compact interval \([s,2s]\). Choose a finite
subcover \(I_{r_1},\ldots,I_{r_N}\), and write \(h_i=h(r_i)\). With
\(\alpha_{r_i}\) defined using this \(h_i\), setting
\[
\beta:=\frac12\min_{1\le i\le N}\alpha_{r_i}>0
\]
gives \eqref{eq:lemzero-uniform-interior}.

Finally, set \(t_0=2s\), \(\delta=\varepsilon/2\), and \(K_0=[\varepsilon/4,M]\). For \(y\in[0,\delta]\),
the strong Markov property at \(\tau_\varepsilon^+\) yields
\[
\begin{split}
P_{t_0}(y,K_0)
&\ge \mathbf E_y\left[
\I_{\{\tau_\varepsilon^+<s,\ X_{\tau_\varepsilon^+}\le M/2\}}
P_{t_0-\tau_\varepsilon^+}(X_{\tau_\varepsilon^+},K_0)\right]\\
&\ge \frac{q\beta}{2}.
\end{split}
\]
Here the last inequality follows from \eqref{eq:lemzero-uniform-interior}, since on the event in the indicator
\(X_{\tau_\varepsilon^+}\in[\varepsilon,M/2]\) and \(t_0-\tau_\varepsilon^+\in[s,2s]\).
Taking \(c_0=q\beta/2\) proves the second assertion.
\qed

\bglemma\label{psiirr}
Suppose that $R_0(0)>0$, $(P_t)_{t\ge0}$ is strong Feller on $(0,\infty)$, and $(P_t)_{t\ge0}$ is irreducible on $(0, \infty)$. Fix \(x_0\in(0,\infty)\) and let
\[
U(x,A):=\int_0^\infty \mrm e^{-t}P_t(x,A)\d t,\qquad
\psi(A):=U(x_0,A),
\]
for \(x\in[0,\infty)\) and \(A\in\mathcal B([0,\infty))\). Then the process is \(\psi\)-irreducible on \([0,\infty)\). Equivalently, if \(\psi(A)>0\), then \(U(x,A)>0\) for all \(x\in[0,\infty)\).
\edlemma
\proof
Suppose \(A \in \mathcal B([0,\infty))\) satisfies \(\psi(A)>0\). We first show that \(U(x,A)>0\) for all \(x\in(0,\infty)\).
For any fixed $A \in \mathcal{B}([0, \infty))$, the strong Feller property implies that the mapping $y \mapsto P_t(y, A)$ is continuous on $(0, \infty)$ for every $t > 0$. By the Dominated Convergence Theorem, the mapping $y \mapsto U(y, A)$ is also continuous on \((0,\infty)\). 
Since $U(x_0, A) = \psi(A) > 0$, the set \(O_A := \{y \in (0, \infty) : U(y, A) > 0\}\) is a nonempty open subset of $(0, \infty)$.

For any $x \in (0, \infty)$, the open-set irreducibility guarantees that $P_t(x, O_A) > 0$ for all $t > 0$. Hence \(U(x,O_A)=\int_0^\infty \mrm e^{-t}P_t(x,O_A)\d t>0\). Then the 2-step transition kernel satisfies
\beqnn
U^2(x, A) \ar=\ar \int_0^\infty U(x, \d y) U(y, A) \ge \int_{O_A} U(x, \d y) U(y, A) > 0.
\eeqnn
On the other hand, applying the Fubini-Tonelli theorem and the Chapman-Kolmogorov equation, we can compute \(U^2(x, A)\) directly:
\beqnn
U^2(x, A) \ar=\ar \int_0^\infty \left( \int_0^\infty \mrm{e}^{-t} P_t(x, \d y) \d t \right) \left( \int_0^\infty \mrm{e}^{-s} P_s(y, A) \d s \right) \cr
\ar=\ar \int_0^\infty \int_0^\infty \mrm{e}^{-(t+s)} P_{t+s}(x, A) \d s \d t \cr
\ar=\ar \int_0^\infty \mrm{e}^{-u} P_u(x, A) \left( \int_0^u \d t \right) \d u \cr
\ar=\ar \int_0^\infty u \mrm{e}^{-u} P_u(x, A) \d u.
\eeqnn
Since \(U^2(x, A) > 0\), the Lebesgue measure of the set $\{u > 0 : P_u(x, A) > 0\}$ must be strictly positive. Consequently, for all \(x\in(0,\infty)\),
\beqnn
U(x, A) = \int_0^\infty \mrm{e}^{-u} P_u(x, A) \d u > 0.
\eeqnn
The boundary point is covered by the same resolvent argument. Indeed, Lemma \ref{lemzero} gives \(P_t(0,O_A)>0\) for every \(t>0\), and hence \(U(0,O_A)>0\). Therefore
\[
U^2(0,A)\ge \int_{O_A}U(0,\d y)U(y,A)>0,
\]
and the preceding argument yields \(U(0,A)>0\). This proves the assertion.
\qed

\bgproposition\label{petite}
Suppose that $R_0(0)>0$, $(P_t)_{t\ge0}$ is strong Feller on $(0,\infty)$, and $(P_t)_{t\ge0}$ is irreducible on $(0, \infty)$. Then any compact set $K\subset [0, \infty)$ is a petite set; that is, there exist a finite measure $a$ on $(0, \infty)$ and a non-trivial measure $\nu$ on $[0, \infty)$ such that 
\[
\int_0^\infty P_t(x, \cdot) a(\d t) \ge \nu(\cdot) \quad \text{for all } x \in K.
\]
\edproposition
\proof
Let \(U\) and \(\psi\) be as in Lemma \ref{psiirr}.

Consider the discrete-time Markov chain defined by the transition kernel \(U(x, \cdot)\). Since the continuous-time process is $\psi$-irreducible, the resolvent chain is also $\psi$-irreducible. According to the general theory of Markov chains (see Meyn and Tweedie \cite[Theorem 5.2.2]{MT93book}), there exists a Borel set $C$ with $\psi(C) > 0$ satisfying a finite-step minorization: for some integer $m \ge 1$, constant $\eta > 0$, and probability measure $\nu_m$,
\[
  U^m(y, \cdot) \ge \eta \nu_m(\cdot),\qquad y \in C.
\]

Let \(K \subset (0, \infty)\) be an arbitrary compact set. Since $\psi(C) > 0$, Lemma \ref{psiirr} gives \(U(x, C) > 0\) for all $x \in (0, \infty)$. 
By the strong Feller property, the mapping \(x \mapsto U(x, C)\) is continuous on $(0, \infty)$. Since $K$ is compact, this strictly positive continuous function attains its minimum on $K$. Thus, there exists a constant \(\epsilon_K > 0\) such that \(U(x, C) \ge \epsilon_K\) for all $x \in K$.

Using the Chapman-Kolmogorov equation for the resolvent chain, we obtain for any $x \in K$ and $A \in \mathcal{B}([0, \infty))$:
\beqnn
U^{m+1}(x, A) \ar=\ar \int_0^\infty U(x, \d y) U^m(y, A) \cr
\ar\ge\ar \int_C U(x, \d y) U^m(y, A) \cr
\ar\ge\ar \int_C U(x, \d y) \eta \nu_m(A) \cr
\ar=\ar \eta U(x, C) \nu_m(A) \ge \eta \epsilon_K \nu_m(A),
\eeqnn
which gives a finite-step minorization for the compact set $K$ with respect to the resolvent chain. In terms of the original continuous-time semigroup, this reads as
\beqlb\label{petitR}
U^{m+1}(x, \cdot)= \int_0^\infty P_t(x, \cdot) a_{m+1}(\d t)\ge \eta \epsilon_K \nu_m(\cdot),\eeqlb
where the sampling distribution $a_{m+1}(\d t)$ is the Gamma$(m+1, 1)$ distribution (i.e., the $(m+1)$-fold convolution of the standard exponential distribution). 
This implies that $K$ is a petite set for the continuous-time process $(X_t)_{t\ge0}$.

We next extend the petite property down to compact sets touching the boundary. By the second assertion of Lemma \ref{lemzero}, choose \(t_0>0\), \(\delta>0\), \(c_0>0\), and a compact interval \(K_0\subset(0,\infty)\) with nonempty interior such that
\[
\inf_{y\in[0,\delta]}P_{t_0}(y,K_0)\ge c_0.
\]
By the preceding paragraph, \(K_0\) is petite and satisfies \eqref{petitR} with \(\epsilon_{K_0}>0\).
Using the Chapman-Kolmogorov equation and \eqref{petitR}, we have for any starting point $y \in [0, \delta]$ and any measurable set $A$:
\beqnn
\int_0^\infty P_{t_0+s}(y, A)a_{m+1}(\d s)
\ar=\ar
\int_0^\infty \left(\int_{[0,\infty)} P_{t_0}(y, \d x)P_s(x, A) \right) a_{m+1}(\d s)\cr
\ar\ge\ar
\int_{K_0} P_{t_0}(y, \d x) \left(\int_0^\infty P_s(x, A) a_{m+1}(\d s)\right)\cr
\ar\ge\ar
\int_{K_0} P_{t_0}(y, \d x) \eta \epsilon_{K_0} \nu_m(A) \cr
\ar\ge\ar
c_0 \eta \epsilon_{K_0} \nu_m(A),
\eeqnn
By introducing a shifted sampling measure $\tilde{a}(\d t) := a_{m+1}(\d t - t_0)\I_{\{t \ge t_0\}}$, the left-hand side is exactly $\int_0^\infty P_t(y, A)\tilde{a}(\d t)$. This yields that the interval $[0, \delta]$ is petite, with minorizing measure proportional to the same \(\nu_m\) used for interior compact sets. Hence, for an arbitrary compact \(K\subset[0,\infty)\), the parts \(K\cap[0,\delta]\) and \(K\cap[\delta,\infty)\) can be combined by taking a convex combination of the corresponding sampling measures and the smaller of the two constants multiplying \(\nu_m\). Thus every compact subset of \([0,\infty)\) is petite.
\qed

\bglemma\label{aperiodic}
Suppose that $R_0(0)>0$, $(P_t)_{t\ge0}$ is strong Feller on $(0,\infty)$, and $(P_t)_{t\ge0}$ is irreducible on $(0, \infty)$. Then the process is aperiodic: there exist a non-trivial measure \(\psi\) for which the process is \(\psi\)-irreducible and a petite set \(C\) with \(\psi(C)>0\) such that
\[
P_t(x,C)>0,\qquad x\in C,\ t>0.
\]
\edlemma
\proof
Let \(U\) and \(\psi\) be as in Lemma \ref{psiirr}. Choose a compact interval \(C=[a,b]\subset(0,\infty)\) with nonempty interior and \(\psi(C)>0\). Then \(C\) is petite by Proposition \ref{petite}.
Moreover, by the open-set irreducibility on \((0,\infty)\), for every \(x\in C\) and every \(t>0\),
\[
P_t(x,C)\ge P_t(x,C^\circ)>0.
\]
This verifies aperiodicity. 
\qed

{\bf Proof of Theorem \ref{thmergo}.}
First, under the stated alternatives \hyperlink{cond:C1}{\textup{(C1)}} or both \hyperlink{cond:C2}{\textup{(C2)}} and \hyperlink{cond:C3}{\textup{(C3)}}, Theorems \ref{SF_a} and \ref{irred} give the interior strong Feller property and irreducibility. Lemma \ref{lemzero} extends the irreducibility input to the boundary point \(0\), using the assumption \(R_0(0)>0\). Hence, Lemma \ref{psiirr} gives \(\psi\)-irreducibility, Proposition \ref{petite} implies that every compact subset of $[0, \infty)$ is a petite set, and Lemma \ref{aperiodic} gives aperiodicity.

To obtain uniform exponential ergodicity, we construct a bounded Lyapunov function to directly verify the continuous-time drift condition in \cite{DMT95}. Define the function $\tilde{V} \in C_b^2([0, \infty))$ by
\beqnn
\tilde{V}(x) = 2 - (x+1)^{-r}, \quad x \ge 0,
\eeqnn
where $r > 0$ is the constant given in Condition \hyperlink{cond:C4}{\textup{(C4)}}. Notice that $1 \le \tilde{V}(x) < 2$ for all $x \ge 0$. The derivatives of $\tilde{V}$ are given by
\beqnn
\tilde{V}'(x) = r(x+1)^{-r-1} > 0, \quad \tilde{V}''(x) = -r(r+1)(x+1)^{-r-2} < 0.
\eeqnn
Applying the  generator $L$ to $\tilde{V}$, we have
\beqlb\label{eq:LV_bound}
L\tilde{V}(x) \ar=\ar R_0(x)\tilde{V}'(x) + \frac{1}{2}R_1(x)\tilde{V}''(x) \cr
\ar\ar + R_2(x) \int_0^\infty \left[ \tilde{V}(x+z) - \tilde{V}(x) - z\tilde{V}'(x) \right] \pi(\d z)\cr
\ar\le\ar
R_0(x)\tilde{V}'(x).
\eeqlb
For the drift term, using Condition \hyperlink{cond:C4}{\textup{(C4)}}, we have
\beqnn
R_0(x)\tilde{V}'(x) \le \left(-\lambda_1 x^{r+1} + \lambda_2\right) r(x+1)^{-r-1}.
\eeqnn
As $x \to \infty$, $R_0(x)\tilde{V}'(x) \to -\lambda_1 r < 0$. Combining all these estimates, we conclude that
\beqnn
\limsup_{x\to\infty} L\tilde{V}(x) \le -\lambda_1 r < 0.
\eeqnn
Since $1 \le \tilde{V}(x) < 2$, we can choose a sufficiently small constant $c \in (0, \lambda_1 r / 2)$ and a sufficiently large $l > 0$ such that for all $x > l$, $L\tilde{V}(x) \le -c \tilde{V}(x)$. For $x \in [0, l]$, since $L\tilde{V}(x)$ is continuous and bounded on compact intervals, there exists a constant $b > 0$ such that for all $x \in [0, \infty)$,
\beqlb\label{eq:drift_final}
L\tilde{V}(x) \le -c\tilde{V}(x) + b\I_{[0, l]}(x).
\eeqlb
Let $K = [0, l]$. Since $K$ is compact, it is a petite set by Proposition \ref{petite}. The inequality \eqref{eq:drift_final} rigorously satisfies the continuous-time drift condition of \cite[Theorem 5.2 (c)]{DMT95}. Therefore, the process $(X_t)_{t\ge0}$ is $\tilde{V}$-uniformly exponentially ergodic. 

Since \(1\le \tilde{V}\le 2\), the $\tilde{V}$-norm is equivalent to the standard total variation norm. It follows that there exist constants $C > 0$ and $\lambda > 0$ such that for any initial distribution $\nu \in \mathcal{P}([0, \infty))$,
\beqnn
\|\nu P_t - \mu^*\|_{\mathrm{TV}} \le C \mrm{e}^{-\lambda t}, \quad t \ge 0,
\eeqnn
where $\mu^*$ is the unique invariant probability measure. This completes the proof.
\qed

\bibliographystyle{abbrv}
\bibliography{Feller527}
\end{document}